\documentclass{cocv}
\usepackage{dsfont}
\usepackage{amsfonts,amsmath,amssymb,amsthm,bm}
\numberwithin{equation}{section}
\usepackage{url}
\begin{document}
\title[Controllability of a $2\times2$  parabolic system]{Controllability of a $2\times2$  parabolic system by one force with space-dependent coupling term of order one.}
\thanks{This work was partially supported by R\'{e}gion de Franche-Comt\'{e} (France).}
%
\author{Michel Duprez}\address{Laboratoire de Math\'{e}matiques de Besan\c{c}on 
UMR CNRS 6623, Universit\'{e} de Franche-Comt\'{e},
16 route de Gray, 25030 Besan\c{c}on Cedex, France, 
E-mail: mduprez@math.cnrs.fr}
\date{January 6, 2016}
\begin{abstract}
 This paper is devoted to the  controllability of  linear systems 
 of two coupled parabolic equations when the coupling involves  a space dependent first order term. 
This system is set on an bounded interval $I \subset \subset \mathbb{R}$, and the first equation is controlled
by a force supported in a subinterval of $I$ or on the boundary. 
In the case where the intersection of the coupling and control domains is nonempty, we prove
null controllability at any time.
Otherwise, we provide a minimal time for null controllability. 
Finally we give  a necessary and sufficient condition for the approximate controllability. 
The main technical tool for  obtaining  these results is  the moment method. 
\end{abstract}
%
%
\subjclass{93B05, 93B07}
\keywords{Controllability, Observability, Moment Method, Parabolic Systems.}
\maketitle
\section{Introduction and main results}
Let $T>0$, $\omega:=(a,b)\subseteq(0,\pi)$ and $Q_T:=(0,\pi)\times(0,T)$. 
We consider in the present paper the  following distributed control system
\begin{equation}\label{system primmal int}
 \left\{\begin{array}{ll}
\partial_ty_1-\partial_{xx}y_1=\mathds{1}_{\omega}v
&\mathrm{in~} Q_T,\\
\partial_ty_2-\partial_{xx}y_2+p(x)\partial_xy_1+q(x)y_1=0
&\mathrm{in~} Q_T,\\
\noalign{\smallskip}y_1(0,\cdot)=y_1(\pi,\cdot)=y_2(0,\cdot)=y_2(\pi,\cdot)=0&\mathrm{on}~(0,T),\\
\noalign{\smallskip}y_1(\cdot,0)=y^0_1,~y_2(\cdot,0)=y^0_2&\mathrm{in}~(0,\pi)
        \end{array}
\right.
\end{equation}
 and boundary control system 
\begin{equation}\label{system primmal bord}
 \left\{\begin{array}{ll}
\partial_tz_1-\partial_{xx}z_1=0&\mathrm{in~} Q_T,\\
\partial_tz_2-\partial_{xx}z_2+p(x)\partial_xz_1+q(x)z_1=0&\mathrm{in~} Q_T,\\
\noalign{\smallskip}z_1(0,\cdot)=u,~z_1(\pi,\cdot)=z_2(0,\cdot)=z_2(\pi,\cdot)=0&\mathrm{on}~(0,T),\\
\noalign{\smallskip}z_1(\cdot,0)=z^0_1,~z_2(\cdot,0)=z^0_2&\mathrm{in}~(0,\pi),
        \end{array}
\right.
\end{equation}
where $y^0:=(y^0_1,y^0_2)\in L^2(0,\pi)^2$ and $z^0:=(z^0_1,z^0_2)\in H^{-1}(0,\pi)^2$ 
are the initial conditions, 
$v\in L^2(Q_T)$ and $u\in L^2(0,T)$ are the controls, $p\in W^1_{\infty}(0,\pi)$,  
 $q\in L^{\infty}(0,\pi)$.

It is known (see \cite[p. 102]{lionscontrole} (resp.  
 \cite[Prop. 2.2]{fernandezcaraboundary2010})) that for given initial data $y^0\in L^2(0,\pi)^2$ (resp. $z^0\in H^{-1}(0,\pi)^2$) 
and a control $v\in L^2(Q_T)$ (resp. $u\in L^2(0,T)$) 
System \eqref{system primmal int} (resp. \eqref{system primmal bord}) has a unique solution 
$y=(y_1,y_2)$ (resp. $z=(z_1,z_2)$) in 
\begin{equation*}\begin{array}{c}
L^2(0,T;H^1_0(0,\pi)^2)\cap \mathcal{C}([0,T];L^2(0,\pi)^2) \vspace*{2mm}\\
(resp.~ L^2(Q_T)^2\cap \mathcal{C}([0,T];H^{-1}(0,\pi)^2)~),
\end{array}\end{equation*}
which depends continuously on the initial data and the control, that is
\begin{equation*}\begin{array}{c}
\|y\|_{L^2(0,T;H^1_0(0,\pi)^2)}+\|y\|_{\mathcal{C}([0,T];L^2(0,\pi)^2)}
\leqslant C_T(\|y^0\|_{L^2(0,\pi)^2}+\|v\|_{L^2(Q_T)})\vspace*{2mm} \\
(resp.~\|z\|_{ L^2(Q_T)^2}+\|z\|_{ \mathcal{C}([0,T];H^{-1}(0,\pi)^2)}
\leqslant C_T(\|z^0\|_{H^{-1}(0,\pi)^2}+\|u\|_{L^2(0,T)})
~),
\end{array}\end{equation*}
where $C_T$ does not depend on $y^0$, $v$, $z^0$ and $u$.
 
Let us introduce the notion of null and approximate controllability for this kind of systems. 
\begin{enumerate}
\item[$\bullet$]
System \eqref{system primmal int} (resp. System \eqref{system primmal bord}) 
is \textit{null controllable} at time $T$ 
if for every initial condition  $y^0\in  L^2(0,\pi)^2$ 
(resp. $z^0\in  H^{-1}(0,\pi)^2$) there exists a control  $v\in L^2(Q_T)$ 
(resp. $u\in L^2(0,T)$) 
such that the solution  to System \eqref{system primmal int} 
(resp. System \eqref{system primmal bord}) satisfies
\begin{equation*}
 y(T)\equiv0 ~~~~~\mathrm{(resp.} ~z(T)\equiv0\mathrm{)~~~~~in~~}(0,\pi).
\end{equation*}
\item[$\bullet$]
System \eqref{system primmal int} (resp. System \eqref{system primmal bord}) is 
\textit{approximately controllable} at time $T$ 
if for all $\varepsilon>0$ and  all   $y^0, ~y^1\in  L^2(0,\pi)^2$ 
(resp. $z^0,~ z^1\in  H^{-1}(0,\pi)^2$)
there exists a control $v\in L^2(Q_T)$ (resp. $u\in L^2(0,T)$)  
such that the solution to System \eqref{system primmal int} 
(resp. System \eqref{system primmal bord}) 
satisfies
\begin{equation*}
 \|y(T)-y^1\|_{L^2(0,\pi)^2}\leqslant \varepsilon
 ~~~~~~\mathrm{(resp. }~ \|z(T)-z^1\|_{H^{-1}(0,\pi)^2}\leqslant \varepsilon\mathrm{)}.
\end{equation*}
\end{enumerate}
We recall  that null-controllability at some time T implies approximate controllability at the same time T for linear parabolic systems.
This follows from the backward uniqueness result of \cite[Th. 1.1]{Ghidaglia86} for first order
perturbations and Propositions 2.5 and 2.6.
Moreover the approximate controllability does not depend on the time of control $T$ 
since we  consider autonomous systems. It is a consequence of the analyticity in time of the adjoint semigroup.

The main goal of this article is to provide a complete answer to the null and approximate controllability issues 
for System \eqref{system primmal int} and \eqref{system primmal bord}. 
For a survey  and some applications in physics, chemistry or biology 
concerning the controllability of this kind of systems, we refer to \cite{ammar2011recent}. 
In the last decade, many papers studied this problem, 
however most of them are related to some parabolic systems with  zero order coupling terms. 
 Without first order coupling terms, some Kalman coupling conditions are made explicit in  \cite[Th. 1.4]{ammar2009generalization}, \cite[Th. 1.1]{ammar2009kalman} and \cite[Th. 1.1]{fernandezcaraboundary2010}
for distributed null controllability of systems of more than two equations with constant matrices and in higher space dimension and  
in the case of time dependent matrices, some Silverman-Meadows coupling conditions 
are given in \cite[Th. 1.2]{ammar2009generalization}.

Concerning the null and approximate controllability of Systems \eqref{system primmal int} and   \eqref{system primmal bord} in the case $p\equiv0$ and $q\not\equiv0$ in $(0,\pi)$,  
a partial answer is given  in 
\cite{alabau2011indirect,alabau_leautaud2013,dehman_leautaud2014,rosier_teresa_2011} under the sign condition 
$q\leqslant0\mathrm{~~or~~}q\geqslant 0~~\mathrm{in~}(0,\pi).$
These results are obtained as a consequence of controllability results of a hyperbolic system using the transmutation method (see \cite{miller_transmutation}). 
One can find a necessary and sufficient condition in \cite{cherifmalaga} when 
$\int_{0}^{\pi}q(x)dx\neq0$.
Finally, in \cite{oliveapprox2014}, the authors gives a complete characterization 
of the approximate controllability and, 
 in the recent work \cite{CherifMinimalTimeDisjoint,Ammar-Khodja2015},
we can find   a complete study of the null controllability. 

When $p\neq0$, the approximate controllability of systems \eqref{system primmal int} and \eqref{system primmal bord} in any dimension is studied in  \cite{olive_bound_appr_2014}. 
 The author gives a  sufficient condition for the approximate controllability 
on the boundary and, in the case of analytic coupling coefficients $p$ and $q$, a necessary and sufficient condition for the internal approximate controllability.

Let us now remind known results concerning null controllability for systems of the following more general form.
Let $\Omega$ be a bounded domain in $\mathbb{R}^N$ ($N\in\mathbb{N}^*$) of class $\mathcal{C}^2$ and $\omega_0$ an arbitrary nonempty subset of $\Omega$. 
We denote by $\partial\Omega$ the boundary of $\Omega$. 
Consider the  system of two coupled linear parabolic equations
\begin{equation}\label{system primmal int gene}
 \left\{\begin{array}{ll}
\partial_ty_1=\Delta y_1+g_{11}\cdot\nabla y_1+g_{12}\cdot\nabla y_2 +a_{11}y_1+a_{12}y_2+\mathds{1}_{\omega_0}v
&\mathrm{in~} \Omega\times(0,T),\\
\noalign{\smallskip}\partial_ty_2=\Delta y_2+g_{21}\cdot\nabla y_1+g_{22}\cdot\nabla y_2 +a_{21}y_1+a_{22}y_2
&\mathrm{in~} \Omega\times(0,T),\\
\noalign{\smallskip}y=0&\mathrm{on}~\partial\Omega\times(0,T),\\
\noalign{\smallskip}y(\cdot,0)=y^0&\mathrm{in}~\Omega,
        \end{array}
\right.
\end{equation}
where $y^0\in L^2(\Omega)^2$, $g_{ij}\in L^{\infty}(\Omega\times(0,T))^N$ and $a_{ij}\in L^{\infty}(\Omega\times(0,T))$ for all $i,j\in\{1,2\}$.

As a particular case of the result in section 4 of \cite{gonzalez2010controllability} (see also \cite{MR2148852}), System \eqref{system primmal int gene} is null controllable whenever 
\begin{equation}\label{cond casc}
 g_{21}\equiv 0\mbox{ in }\omega_1~~\mathrm{and}~~
 (a_{21}>C \mbox{ in }\omega_1\mathrm{~or~}a_{21}<-C\mbox{ in }\omega_1) ,
\end{equation}
for a  positive constant $C$ and $\omega_1$ a non-empty open subset of $\omega_0$.

In \cite[Th. 4]{guerrerosyst22}, the author supposes that $a_{11},~g_{11},~a_{22},~g_{22}$ are constant 
and the first order coupling operator $g_{21}\cdot\nabla+a_{21}$ can be written as
\begin{equation}\label{form guerrero}
g_{21}\cdot\nabla+a_{21}=P_1\circ\theta \mathrm{~~in~~}\Omega\times(0,T),
\end{equation}
where $\theta\in \mathcal{C}^2(\overline{\Omega})$ satisfies $|\theta|>C$ in $\omega_1\subseteq\omega_0$ for a  positive constant $C$ and $P_1$ is given by 
$P_1:=m_0\cdot\nabla+m_1$,
for some $m_0,m_1\in\mathbb{R}$. 
Moreover the operator $P_1$ has to satisfy
\begin{equation*}
 \|u\|_{H^1(\Omega)}\leqslant C \|P^*_1u\|_{L^2(\Omega)}~~~~\forall u\in H^1_0(\Omega).
\end{equation*}
Under these assumptions, the author proves the null controllability of System \eqref{system primmal int gene} at any time.

In \cite[Th. 2.1]{benabdallah2014}, the authors prove that the same property holds true for System \eqref{system primmal int gene} if we assume that 
$a_{ij}\in \mathcal{C}^4(\overline{\Omega\times(0,T)})$, $g_{ij}\in \mathcal{C}^1(\overline{\Omega\times(0,T)})^N$ for all $i,j\in\{1,2\}$, $g_{21}\in \mathcal{C}^3(\overline{\Omega\times(0,T)})$  and 
 the geometrical condition 
\begin{equation}\label{cond assia}
\left\{\begin{array}{l}
\partial\omega\cap\partial\Omega\mbox{ contains an open subset }\gamma
\mbox{ for which the interior }\mathring{\gamma}\mbox{ is non-empty},\\
 \noalign{\smallskip}\exists x_0\in\gamma \mathrm{~s.t.~}g_{21}(t,x_0)\cdot \nu(x_0)\neq0 \mathrm{~~for~all~}t\in[0,T],
\end{array}\right.
\end{equation}
where $\nu$ represents the exterior normal unit vector  to the boundary $\partial\Omega$.

Lastly, for constant coefficients, it is proved  in 
\cite[Th. 1]{pierreordre1} 
that System \eqref{system primmal int gene} is null/approximately controllable at any time $T$ 
if and only if 
\begin{equation*}
 g_{21}\neq0\mathrm{~~~or~~~}a_{21}\neq0.
\end{equation*}
In \cite[Th. 2]{pierreordre1}, the authors give also a condition of null/approximate controllability in dimension one 
which can be written for system \eqref{system primmal int} as: $p\in\mathcal{C}^2(\omega_0)$, $q\in\mathcal{C}^3(\omega_0)$ and 
   \begin{equation*}
 \begin{array}{c}  -4 \partial_x(q) \partial_x(p)p+  \partial_{xx}(q)p^2+2q\partial_x(q)p
      -3pq \partial_{xx}p +6 q   (\partial_xp)^2-2 q^2\partial_xp\\
     -  \partial_{xxx}(p)p^2 +5    \partial_x(p)\partial_{xx}(p)p  -  4(\partial_xp)^3\neq0\mbox{ in }\omega_0
     \end{array}
  \end{equation*} 
for a subinterval $\omega_0$ of $\omega$. 

Now let us go back to Systems  \eqref{system primmal int} and \eqref{system primmal bord} for which we will provide a  complete description of the null and approximate controllability. 
Our first and main result is the following
\begin{thrm}\label{theo null int}
Let us suppose that $p\in W^1_{\infty}(0,\pi)\cap W^2_{\infty}(\omega)$, 
$q\in L^{\infty}(0,\pi)\cap W^1_{\infty}(\omega)$ and 
\begin{equation}\label{Ik Ik1 neq 0}
\begin{array}{c}
 (\mathrm{Supp}~p\cup\mathrm{Supp}~q)\cap\omega\neq\varnothing.
\end{array}
\end{equation}
Then System \eqref{system primmal int} is null controllable at any time $T$.
\end{thrm}

Let us compare this result with the previously described results to highlight our main contribution:
\begin{enumerate}
 \item Even though System \eqref{system primmal int} is considered in one  space dimension, we remark first that our coupling operator has a more general form than the one in \eqref{form guerrero}  
  assumed in \cite{guerrerosyst22}. 
  Moreover unlike \cite{pierreordre1}, its coefficients are non-constant with respect to the space variable. 
\item We do not have 
the geometrical restriction \eqref{cond assia} assumed in \cite{benabdallah2014}.
More precisely we do not require the control  support to be a neighbourhood of a part of the boundary.
\item  As said before, in \cite[Theorem 4.1]{olive_bound_appr_2014}, the author gives a necessary and sufficient condition for the approximate controllability of System 
\eqref{system primmal int} when $p$ and $q$ are analytic. 
We deduce that the condition of \cite{olive_bound_appr_2014} is satisfied in dimension  one 
as soon as $p$ or $q$ is not equal to zero. 
\end{enumerate}


For all  $k\in\mathbb{N}^*$, 
we denote by  $\varphi_k:x\mapsto \sqrt{\frac{2}{\pi}}\sin(kx)$ 
  the normalized eigenvector of the 
Laplacian operator, with Dirichlet boundary condition, 
and   consider the two following quantities
\begin{equation}\label{def Ik}
 I_{a,k}(p,q):=\displaystyle\int_0^{a}
 \textstyle\left(q-\frac12\partial_xp\right)\varphi_k^2
\hspace*{5mm} \mbox{and}\hspace*{5mm} 
I_k(p,q):=\displaystyle\int_0^{\pi} \textstyle\left(q-\frac12\partial_xp\right)\varphi_k^2,
\end{equation}
for all $k\in \mathbb{N}^*$. Combined with the criterion of Fattorini (see \cite[Cor. 3.3]{fattorini66} or Th. \ref{theo hautus} in the present paper),  Theorem \ref{theo null int}
leads to the following characterization:
\begin{thrm}\label{theo approx int}
Let us suppose that $p\in W^1_{\infty}(0,\pi)\cap W^2_{\infty}(\omega)$ and 
$q\in L^{\infty}(0,\pi)\cap W^1_{\infty}(\omega)$. 
System \eqref{system primmal int} is approximately controllable at any time $T>0$ 
if and only if 
\begin{equation}\label{cond approx1}
\begin{array}{c}
(\mathrm{Supp}~p\cup\mathrm{Supp}~q)\cap\omega\neq\varnothing
\end{array}
\end{equation}
\begin{equation*}
or
\end{equation*}
\begin{equation}\label{cond approx2}
\begin{array}{c}
  |I_k(p,q)|+|I_{a,k}(p,q)|\neq0\mathrm{~for~all~}k\in\mathbb{N}^*.
\end{array}
\end{equation}
\end{thrm}

This last result recovers the case  $p\equiv0$ 
 studied in \cite{oliveapprox2014} for $\mathrm{Supp}~q\cap \omega=\varnothing$, where the authors also use  the criterion of Fattorini.

 \begin{rmrk}
We will see in the proof of Theorems \ref{theo null int} and \ref{theo approx int} that only the following regularity are needed for $p$ and $q$
\begin{equation*}\left\{\begin{array}{l}
 p\in W^1_{\infty}(0,\pi)\cap W^2_{\infty}(\widetilde{\omega}),\\\noalign{\smallskip} 
q\in L^{\infty}(0,\pi)\cap W^1_{\infty}(\widetilde{\omega}),
\end{array}\right.\end{equation*}
for an open subinterval $\widetilde{\omega}$ of $\omega$. These hypotheses are used in Definition \eqref{def Ik} of $I_k(p,q)$ and $I_{a,k}(p,q)$ 
and the change of unknown described in Section \ref{sec: resol}. For more general coupling terms, these control problems are open.
 \end{rmrk}

When the supports of the control and the coupling terms are disjoint in System \eqref{system primmal int}, following the ideas in \cite[Th. 1.3]{Ammar-Khodja2015} 
where the authors studied the case $p\equiv0$, 
we obtain a minimal time of null controllability:
\begin{thrm}\label{theo null T0}
Let $p\in W^1_{\infty}(0,\pi)$, 
$q\in L^{\infty}(0,\pi)$. Suppose that Condition  \eqref{cond approx2} holds and
\begin{equation}\label{cond p q non nuls}
 (\mathrm{Supp}~p\cup\mathrm{Supp}~q)\cap\omega=\varnothing.
 \end{equation}
Let $T_0(p,q)$ be given by 
\begin{equation}\label{def T_0}
 T_0(p,q):=\limsup\limits_{\substack{k \to \infty }} 
 \frac{\min(-\log \left| I_k(p,q) \right|,-\log \left| I_{a,k}(p,q) \right|) }{k^2} .
\end{equation}
One has
\begin{enumerate}
 \item If $T>T_0(p,q)$, then System \eqref{system primmal int} 
 is null controllable at time $T$. 
 \item If $T<T_0(p,q)$, then System \eqref{system primmal int} is not null controllable at time $T$. 
\end{enumerate}
\end{thrm}

Concerning the boundary controllability, in \cite[Th. 3.3]{olive_bound_appr_2014}, 
using the criterion of Fattorini, the author proves that 
System \eqref{system primmal bord} is approximately controllable at time $T$ if and only if 
  \begin{equation}\label{Ik neq 0}
   I_k(p,q)\neq0\mathrm{~for~all~}k\in\mathbb{N}^*.
  \end{equation}
About null controllability of System \eqref{system primmal bord}, we can again generalize the results given in \cite[Th. 1.1]{Ammar-Khodja2015} to obtain a minimal time:
\begin{thrm}\label{theo bord}
Let $p\in W^1_{\infty}(0,\pi)$,  $q\in L^{\infty}(0,\pi)$ and suppose that 
Condition \eqref{Ik neq 0} is satisfied. Let us define
\begin{equation}\label{def T1234}
 T_1(p,q):= \limsup\limits_{\substack{k \to \infty }} 
 \frac{-\log \left| I_k(p,q) \right| }{k^2} .
\end{equation}
 One has 
 \begin{enumerate}
\item If $T>T_1(p,q)$, 
then System \eqref{system primmal bord} is null controllable at time $T$. 
\item If $T<T_1(p,q)$, 
then System \eqref{system primmal bord} is not null controllable at time $T$. 
 \end{enumerate}

\end{thrm}

\begin{rmrk}Using Riemann-Lebesgue Lemma,   sequences 
$(I_k(p,q))_{k\in\mathbb{N}^*}$ and $(I_{a,k}(p,q))_{k\in\mathbb{N}^*}$ are convergent, 
more precisely 
\begin{equation*}
\lim\limits_{k\rightarrow\infty}I_k(p,q)=I(p,q):=\dfrac{1}{\pi}\displaystyle\int_0^{\pi}(q-\textstyle\frac12\partial_xp)
\hspace*{5mm} \mbox{and}\hspace*{5mm} 
\lim\limits_{k\rightarrow\infty}I_{a,k}(p,q)=I_a(p,q):=\dfrac{1}{\pi}\displaystyle\int_0^{a}(q-\textstyle\frac12\partial_xp).
\end{equation*}
Thus, if one of the two 
 limits $I(p,q)$ and $I_a(p,q)$ (resp. the first limit) is not equal  to zero, 
 then the minimal time $T_0(p,q)$ (resp. $T_1(p,q)$) is equal to zero.

\end{rmrk}

This article is organized as follows. In the second section, we present some preliminary results useful 
to reduce the null controllability issues to  the moment problem. In the third and fourth sections, we study the null controllability issue of System \eqref{system primmal int} in the two cases when the intersection of the coupling and control supports is empty or not. 
Then we give the proof of Theorems \ref{theo approx int} and \ref{theo bord} in Section 5 and 6, 
respectively.  
We finish with some comments and open problems in Section 7.

\section{Preliminary results}
Consider the differential operator
\begin{equation*}
\begin{array}{crcl}
 L:&D(L)\subset L^2(0,\pi)^2&\rightarrow &L^2(0,\pi)^2\\
 &f&\mapsto& -\partial_{xx}f+A_0(p\partial_xf+qf),
\end{array}
\end{equation*}
where the matrix $A_0$ is given by
\begin{equation*}
 A_0:=\left(\begin{array}{cc}0&0\\1&0\end{array}\right),
\end{equation*}
the domain of $L$ and its adjoint $L^*$ is given by $D(L)=D(L^*)=H^2(0,\pi)^2\cap H^1_0(0,\pi)^2$. 
In section \ref{sec:biortho}, we will first establish some properties of  the operator $L$ that will be useful for the moment method and, in section \ref{sec:dual}, 
we will recall some characterizations 
of the approximate and null controllability of system \eqref{system primmal int}.

\subsection{Biorthogonal basis}\label{sec:biortho}

Let us first analyze the spectrum of the operators $L$ and $L^*$.
\begin{prpstn}\label{prop base}
For all $k\in\mathbb{N}^*$  consider the two vectors
 \begin{equation*}
\Phi_{1,k}^*:=\left(\begin{array}{c}\psi_k^*\\\varphi_k\end{array}\right),
\Phi_{2,k}^*:=\left(\begin{array}{c}\varphi_k\\0\end{array}\right),
\end{equation*}
where $\psi_k^*$ is defined for all $x\in (0,\pi)$ by
  \begin{equation*}\left\{\begin{array}{l}
 \psi_k^*(x)=\alpha_k^*\varphi_k(x)-\dfrac{1}{k}\displaystyle\int_0^x\sin(k(x-\xi))
 [I_k(p,q)\varphi_k(\xi)+\partial_x(p(\xi)\varphi_k(\xi))-q(\xi)\varphi_k(\xi)]d\xi,\\\noalign{\smallskip}
 \alpha_k^*=\dfrac{1}{k}\displaystyle\int_0^{\pi}\displaystyle\int_0^x\sin(k(x-\xi))
 [I_k(p,q)\varphi_k(\xi)+\partial_x(p(\xi)\varphi_k(\xi))-q(\xi)\varphi_k(\xi)]\varphi_k(x)d\xi dx.
\end{array}\right.\end{equation*} 
 One has 
  \begin{enumerate} \item The spectrum of  $L^*$ is given by $\sigma(L^*)=\{k^2:k\in\mathbb{N}^*\}$. 
  \item For $k\geqslant 1$, the eigenvalue $k^2$ of $L^*$ is simple 
  (algebraic multiplicity 1) if and only if $I_k(p,q)\neq0$. 
 In this case,  $\Phi_{2,k}^*$ and $\Phi_{1,k}^*$ are  respectively  an eigenfunction 
and a generalized eigenfunction of the operator  $L^*$ associated with the eigenvalue $k^2$, more precisely
\begin{equation}\label{equa vec prop}
\left\{\begin{array}{l}
 (L^*-k^2Id)\Phi_{1,k}^*=I_k\Phi_{2,k}^*,\\\noalign{\smallskip}
  (L^*-k^2Id)\Phi_{2,k}^*=0.
\end{array}\right.
\end{equation}
\item For $k\geqslant 1$, the eigenvalue $k^2$ of $L^*$ is double (algebraic multiplicity 2) 
if and only if $I_k(p,q)=0$. 
In this case,  $\Phi_{1,k}^*$ and $\Phi_{2,k}^*$ are two eigenfunctions 
 of the operator  $L^*$ associated with the eigenvalue $k^2$, that is for $i=1,2$
 \begin{equation*}
  (L^*-k^2Id)\Phi_{i,k}^*=0.
\end{equation*}
 
 \end{enumerate}
\end{prpstn}

\begin{proof}The adjoint operator $L^*$ of $L$ is given by $D(L^*)=D(L)$ 
and $L^*f=-\partial_{xx}f+A_0^*(-\partial_x(pf)+qf)$. 
We can remark first that the resolvent of $L^*$ is compact. Thus the spectrum of $L^*$ reduces to its point spectrum. 
The  eigenvalue problem associated with the operator $L^*$ is
\begin{equation}\label{syst 1 preuve base}
\left\{ 
\begin{array}{ll}
- \partial_{xx}\psi -\partial_x(p(x)\varphi)+ q(x) \varphi = \lambda \psi &\hbox{ in }  ( 0,\pi ) , \\ \noalign{\smallskip}
- \partial_{xx}\varphi = \lambda \varphi
& \hbox{ in }  ( 0,\pi ) , \\ \noalign{\smallskip}
\varphi (0) =\psi (0)= \varphi (\pi ) =\psi (\pi )= 0,&
\end{array}
\right.
\end{equation}
where $(\psi,\varphi)\in D(L^*)$ and $\lambda\in\mathbb{C}$. 
 For $\varphi\equiv0$ in $(0,\pi)$ and $\psi=\varphi_k$ in $(0,\pi)$, 
  $\lambda=k^2$ is an eigenvalue of $L^*$  and 
the vector $\Phi_{2,k}^*:=(\varphi_k,0)$ is an associated  eigenfunction. 
If now $\varphi\not\equiv0$ in $(0,\pi)$, then $\lambda=k^2$ is an eigenvalue 
and $\varphi=\kappa\varphi_k$ with $\kappa\in\mathbb{R}^*$. 
We remark that 
 System \eqref{syst 1 preuve base} has a solution if and only if $I_k(p,q)=0$. 
 If $I_k(p,q)=0$, $\Phi_{1,k}^*:=(\psi_k^*,\varphi_k)$ is a second eigenfunction of $L^*$ linearly independent of 
 $\Phi_{2,k}^*$, where,
applying the Fredholm alternative,  
  $\psi_k^*$ is the unique solution to the non-homogeneous Sturm-Liouville problem
\begin{equation}\label{sturm liouville}
\left\{ \begin{array}{l}
-\partial_{xx}\psi - k^2 \psi = f\hbox{ in }  ( 0,\pi ) , \\ 
\noalign{\smallskip}
\psi (0 )=\psi (\pi )=0 , ~ \int_{0}^{\pi }\psi (x) \varphi _{k}(x) \, dx=0. \\ 
\end{array}\right. 
\end{equation}
with
  $ f:= \partial_x(p(x)\varphi_k)- q(x)  \varphi _{k}$. 
We recall that $\int_0^{\pi}f\varphi_k=0$, since $I_k=0$.
Solving System \eqref{sturm liouville} 
leads to the expression of $\psi_k^*$ 
given in  Proposition \ref{prop base}.  
The expression of $\alpha_k$ is given by the equality 
$\int_{0}^{\pi }\psi^*_k (x) \varphi _{k}(x) \, dx=0$
and the identity $\int_0^{\pi}f\varphi_k=0$ leads to $\psi_k(\pi)=0$. 
 Thus, in the case $I_k(p,q)=0$, $\lambda=k^2$ is a double eigenvalue of $L^*$. 
 Items 1 and 3 are now proved.
 
 Let us now suppose that $I_k(p,q)\neq0$. The eigenvalue $\lambda=k^2$ 
 is simple, $\Phi_{2,k}^*:=(\varphi_k,0)$ is an eigenfunction and 
 a solution $\Phi_{1,k}^*:=(\psi,\varphi)$ to $(L^*-k^2Id)\Phi_{1,k}^*=I_k(p,q)\Phi_{2,k}^*$, that is 
 \begin{equation}\label{syst 2 preuve base}
\left\{ 
\begin{array}{ll}
- \partial_{xx}\psi -\partial_x(p(x)\varphi)+ q(x) \varphi = k^2 \psi+I_k(p,q)\varphi_k &\hbox{ in }  ( 0,\pi ) , \\ \noalign{\smallskip}
- \partial_{xx}\varphi = k^2 \varphi& \hbox{ in }  ( 0,\pi ) , \\ \noalign{\smallskip}
\varphi (0) =\psi (0)= \varphi (\pi ) =\psi (\pi )= 0,&
\end{array}
\right.
\end{equation}
 is a generalized eigenfunction of $L^*$. 
  We deduce that $\varphi=\varphi_k$ in $(0,\pi)$ 
 and $\psi$ is solution to the Sturm-Liouville problem \eqref{sturm liouville} with
 $ f=I_k(p,q)\varphi_k+\partial_x(p(x)\varphi_k)- q(x) \varphi_k $.
We obtain the expression of $\psi_k^*$ given 
in  Proposition \ref{prop base}.  
\end{proof}

The function $\psi_k^*$ given in Proposition \ref{prop base} will play an important role in this paper. 
Since $\varphi_k$, $\varphi_k'$
and $I_k$ are bounded, we have the following lemma
\begin{lmm} There exists a positive constant $C$ such that
 \begin{equation}\label{estm psi2}
 |\alpha_k^*|\leqslant \dfrac{C}{k},\hspace*{5mm} 
 \|\psi_k^*\|_{L^{\infty}(0,\pi)}\leqslant \dfrac{C}{k},\hspace*{5mm} 
 \|\partial_x\psi_k^*\|_{L^{\infty}(0,\pi)}\leqslant C,\hspace*{5mm} 
 \forall k\in\mathbb{N}^*.
 \end{equation}
\end{lmm}

 Since the eigenvalues of the operator $L^*$ are real, we deduce that $L$ and $L^*$ 
 have the same spectrum and   
 the associated eigenspaces have the same dimension. The eigenfunctions and the generalized eigenfunctions of $L$ 
 can be found as previously.
 \begin{prpstn}\label{prop base0}
For all $k\in\mathbb{N}^*$  consider the two vectors
\begin{equation*}
\Phi_{1,k}:=\left(\begin{array}{c}0\\\varphi_k\end{array}\right),
\Phi_{2,k}:=\left(\begin{array}{c}\varphi_k\\\psi_k\end{array}\right),
\end{equation*}
  where $\psi_k$ is defined for all $x\in (0,\pi)$ by
  \begin{equation*}\left\{\begin{array}{l}
 \psi_k(x):=\alpha_k\varphi_k(x)-\dfrac{1}{k}\displaystyle\int_0^x\sin(k(x-\xi))
 [I_k(p,q)\varphi_k(\xi)-p(\xi)\partial_x(\varphi_k(\xi))-q(\xi)\varphi_k(\xi)]d\xi,\\  \noalign{\smallskip}
 \alpha_k:=\dfrac{1}{k}\displaystyle\int_0^{\pi}\displaystyle\int_0^x\sin(k(x-\xi))
 [I_k(p,q)\varphi_k(\xi)-p(\xi)\partial_x(\varphi_k(\xi))-q(\xi)\varphi_k(\xi)]\varphi_k(x)d\xi dx,
\end{array}\right.\end{equation*} 
 One has 
  \begin{enumerate}
  \item The spectrum of $L$  is given by $\sigma(L)=\sigma(L^*)=\{k^2:k\in\mathbb{N}^*\}$. 
  \item For $k\geqslant 1$, the eigenvalue $k^2$ of $L$ is simple (algebraic multiplicity 1) if and only if $I_k(p,q)\neq0$. 
 In this case,  $\Phi_{1,k}$ and $\Phi_{2,k}$ are an eigenfunction 
and a generalized eigenfunction of the operator  $L$ associated with the eigenvalue $k^2$, more precisely
\begin{equation}\label{equa vec prop0}
\left\{\begin{array}{l} 
(L-k^2Id)\Phi_{1,k}=0,\\ \noalign{\smallskip}
 (L-k^2Id)\Phi_{2,k}=I_k\Phi_{1,k}.
\end{array}\right.
\end{equation}
\item For $k\geqslant 1$, the eigenvalue $k^2$ of $L$ is double (algebraic multiplicity 2) 
if and only if $I_k(p,q)=0$. 
In this case,  $\Phi_{1,k}$ and $\Phi_{2,k}$ are two eigenfunctions 
 of the operator  $L$ associated with the eigenvalue $k^2$, that is for $i=1,2$
 \begin{equation*}
  (L-k^2Id)\Phi_{i,k}=0.
\end{equation*}
 \end{enumerate}
\end{prpstn}

Lemma 2.3 and Corollary 2.6 in \cite{Ammar-Khodja2015} can be adapted 
easily to prove the following proposition.

\begin{prpstn}
 Consider the families
\begin{equation*}
\begin{array}{l}\mathcal{B}:=\left\{\Phi_{1,k},\Phi_{2,k}:k\in\mathbb{N}^*\right\}
\hspace*{5mm} \mathrm{and}\hspace*{5mm} 
\mathcal{B}^*:=\left\{\Phi^*_{1,k},\Phi^*_{2,k}:k\in\mathbb{N}^*\right\}.
 \end{array}
\end{equation*}
Then
\begin{enumerate}
 \item The sequences $\mathcal{B}$ and $\mathcal{B}^*$ are biorthogonal Riesz bases of $L^2(0,\pi)^2$.
   \item The sequence $\mathcal{B}^*$ is a Schauder basis  of $H^1_0(0,\pi)^2$ and $\mathcal{B}$ is its biorthogonal basis in $H^{-1}(0,\pi)$.
\end{enumerate}
\end{prpstn}

We recall that $\mathcal{B}$ and $\mathcal{B}^*$ are biorthogonal in $L^2(0,\pi)^2$ if
$\langle \Phi_{i,k},\Phi_{j,l}^*  \rangle_{L^2(0,\pi)^2}=\delta_{i,j}\delta_{k,l}$
 for all $k,l\in\mathbb{N}^*$ and $i,j\in\{1,2\}$.

\subsection{Duality}\label{sec:dual}

As it is well known, the controllability 
has a dual concept called \textit{observability} (see for
instance \cite[Th. 2.1]{ammar2011recent}, \cite[Th. 2.44, p. 56–57]{coron2009control}). 
Consider the dual system associated with  System \eqref{system primmal int}
\begin{equation}\label{system dual}
 \left\{\begin{array}{ll}
-\partial_t\theta-\partial_{xx}\theta+A_0^*(-\partial_x(p(x)\theta)+q(x)\theta)=0&\mathrm{in}~ Q_T,\\  \noalign{\smallskip}
\theta(0,\cdot)=\theta(\pi,\cdot)=0&\mathrm{on}~(0,T),\\  \noalign{\smallskip}
\theta(\cdot,T)=\theta^0&\mathrm{in}~(0,\pi),
        \end{array}
\right.
\end{equation}
where $\theta^0\in L^2(0,\pi)^2$. 
Let $B$ the matrix given by $B=(1~ 0)^*$. 
The approximate controllability is equivalent to a \textit{unique continuation property}:

\begin{prpstn}
\begin{enumerate}
 \item System \eqref{system primmal int} is approximately controllable at time $T$
 if and only if 
 for all initial condition $\theta^0\in L^2(0,\pi)^2$ the solution to System 
 \eqref{system dual} satisfies the unique continuation property
 \begin{equation}\label{cont uni int}
  \mathds{1}_{\omega}B^*\theta\equiv0\mathrm{~in~}Q_T~\Rightarrow~\theta\equiv0~\mathrm{in}~Q_T.
 \end{equation}
 \item   System \eqref{system primmal bord} is approximately controllable at time $T$
 if and only if 
 for all initial condition $\theta^0\in H_0^1(0,\pi)^2$ 
 the solution to System \eqref{system dual} satisfies the unique continuation property
 \begin{equation}\label{cont uni bord}
 B^*\partial_x\theta(0,t)\equiv0\mathrm{~in~}(0,T)~\Rightarrow~\theta\equiv0~\mathrm{in}~Q_T.
 \end{equation}
\end{enumerate}
\end{prpstn}

\noindent The null controllability is characterized by an \textit{observability inequality}:

\begin{prpstn}\label{prop ine obs}
\begin{enumerate}
 \item System \eqref{system primmal int} is null controllable at time $T$ if and only if 
 there exists a constant $C_{obs}$ such that 
 for all initial condition $\theta^0\in L^2(0,\pi)^2$ the solution to System \eqref{system dual} 
 satisfies the observability inequality
 \begin{equation}\label{null cont int}
 \|\theta(0)\|^2_{L^2(0,\pi)^2}\leqslant C_{obs} \iint_{Q_T}|\mathds{1}_{\omega}(x)B^*\theta(x,t)|^2dxdt.
 \end{equation}
 \item System \eqref{system primmal int} is null controllable at time $T$ if and only if 
 there exists a constant $C_{obs}$ such that 
 for all initial condition $\theta^0\in H^1_0(0,\pi)^2$ the solution to System \eqref{system dual} 
 satisfies the observability inequality
 \begin{equation}\label{null cont bord}
 \|\theta(0)\|^2_{H^1_0(0,\pi)^2}\leqslant C_{obs} \int_{0}^T|B^*\partial_x\theta(0,t)|^2dt.
 \end{equation} 
\end{enumerate}
\end{prpstn}

\section{Resolution of the moment problem}\label{sec: proof th 1.1}
 In this section, we first establish the moment problem related to the null controllability for System \eqref{system primmal int} and then  we will solve it 
 in section \ref{sec: resol} (Theorem 1.1). 
 The strategy involves finding an equivalent system (see Definition \ref{def:equiv}) to System \eqref{system primmal int},  which has an associated quantity $I_k$ satisfying "some  good properties".

\subsection{Reduction to a moment problem}\label{moment problem}

Let $y^0:=(y_1^0,y_2^0)\in L^2(0,\pi)^2$. For  $i\in \{1,2\}$ and 
 $k\in\mathbb{N}^*$, if we consider $\theta^0:= \Phi_{i,k}^*$ 
in  the dual System \eqref{system dual}, we get after an integration by parts
\begin{equation*} 
\iint_{Q_{T}} v(x,t) \mathds{1}_\omega(x) B^*\theta (x,t)dx dt 
=\langle y(T) , \Phi_{i,k}^* \rangle_{L^2(0,\pi)^2} - \langle y^0 , \theta (0) \rangle_{L^2(0,\pi)^2}.
\end{equation*}
Since $\mathcal{B}^*$ is a Riesz basis of $L^2(0,\pi)^2$, 
System \eqref{system primmal int} is null controllable if and only if 
for all $y^0\in L^2(0,\pi)^2$, there exists a control $v\in L^2(Q_T)$ 
such that for all $k\in\mathbb{N}^*$ and $i\in\{1,2\}$ 
the solution $y$ to System \eqref{system primmal int} satisfies the following equality
\begin{equation} \label{moment}
\iint_{Q_{T}} v(x,t) \mathds{1}_\omega(x) B^{\ast }\theta_{i,k} (x,t) \,dx \, dt 
= - \langle y^0 , \theta_{i,k} (0) \rangle_{L^2(0,\pi)^2},
\end{equation}
where $\theta_{i,k}$ is the solution to the dual system \eqref{system dual} 
with the initial condition $\theta^0 := \Phi_{i,k}^*$.

In the moment problem \eqref{moment}, we will look  for a control $v$ of the form
    \begin{equation}\label{contv} 
v(x,t):=f^{(1)}(x)v^{(1)}(T - t)+f^{(2)}(x)v^{(2)}(T - t) \mathrm{~for~all~} (x,t) \in Q_T,
    \end{equation}
with $v^{(1)}, v^{(2)} \in L^2(0,T)$ and $f^{(1)}, f^{(2)} \in L^2 (0, \pi) $  satisfying 
$\mathrm{Supp} ~f^{(1)}\subseteq \omega$ and $\mathrm{Supp} ~f^{(2)} \subseteq \omega$.  

The solutions $\theta_{1,k}$ and $\theta_{2,k}$ 
 to the dual System \eqref{system dual} with the initial condition  
 $ \Phi_{1,k}^*$ and $ \Phi_{2,k}^*$ 
 are given for all $(x,t)\in Q_T$ by 
\begin{equation}\label{expr thera k 1 2}\left\{\begin{array}{l}
\theta_{1,k}(x,t)=e^{-k^2 (T-t)}\left(\Phi_{1,k}^*(x)-(T-t)I_k(p,q)\Phi_{2,k}^*(x)\right),\\  \noalign{\smallskip}
\theta_{2,k} (x, t) = e^{-k^2 (T-t)}\Phi_{2,k}^*(x). 
\end{array}\right.\end{equation}
Plugging \eqref{contv} and  \eqref{expr thera k 1 2}  
 in  the moment problem \eqref{moment}, we get for all $k\geqslant1$
\begin{equation*}
\left\{\begin{array}{l}
\displaystyle  \widetilde f_k^{(1)}\int_0^Tv^{(1)}(t) e^{- k^2 t}\,dt
+ \widetilde f_k^{(2)}\int_0^Tv^{(2)}(t) e^{- k^2 t}\,dt \\
\noalign{\smallskip}
\displaystyle \hspace*{3cm}
- I_k (p,q) f_k^{(1)}\int_0^T v^{(1)}(t) t e^{-k^2 t} \,dt - I_k(p,q) f_k^{(2)}\int_0^Tv^{(2)}(t) t e^{-k^2 t}\,dt  \\
\noalign{\smallskip}
\displaystyle \hspace*{6cm}
= -e^{-k^2 T}\left\{    y_{1,k}^0 - T I_k (p,q) y_{2,k}^0\right\},\\
 f_k^{(1)}\displaystyle\int_0^T v^{(1)}(t) e^{- k^2 t}\,dt + f_k^{(2)}\int_0^T v^{(2)}(t) e^{- k^2 t}\,dt 
= - e^{-k^2 T} y_{2,k}^0,
\end{array}\right.
\end{equation*} 
where  $f_{k}^{(i)}$, $\widetilde{f}_{k}^{(i)}$ and ${y}_{i,k}^0$ 
are given for all $i\in\{1,2\}$ and $k\in\mathbb{N}^*$ by
\begin{equation}\label{def fik}
f_{k}^{(i)}:=\int_0^\pi  f^{(i)}(x) \varphi_k(x)dx,~~~
\widetilde f_{k}^{(i)}:= \int_0^\pi f^{(i)}(x) \psi_k^* (x)dx,
\end{equation}
and
\begin{equation}\label{def y0ik}
y_{i,k}^0:=\langle y^0,\Phi^*_{i,k}\rangle_{L^2(0,\pi)}.
\end{equation}

In \cite[Prop. 4.1]{fernandezcaraboundary2010}, the authors proved that 
the family $\left\{ e_{1,k}:=e^{-k^{2}t}, e_{2,k} := te^{-k^{2}t} \right\}_{k\geq 1} $ admits a biorthogonal 
family $\{q_{1,k},q_{2,k}\}_{k\geq 1}$ in the space $L^{2}(0,T)$,  \emph{i.e.} a family satisfying
\begin{equation}\label{orth bi base}
\int_{0}^{T}e_{i,k}q_{j,l}(t) \, dt = \delta_{ij}\delta_{kl}, \quad \forall  k,l\geq 1, \quad 1\leq i,j\leq 2. 
\end{equation}
Moreover  for all $\varepsilon >0$ there exists a constant $C_{\varepsilon ,T}>0$ such that 
\begin{equation}\label{estim qik}
\| q_{i,k} \Vert_{L^{2}( 0,T) }\leq C_{\varepsilon ,T} \, e^{\varepsilon k^2}, \quad \forall k \ge 1, \quad i=1,2.
\end{equation}
We will look for $v^{(1)}$ and $v^{(2)}$ of the form 
\begin{equation}\label{contvi} 
 v^{(i)}(t)=\sum\limits_{k\geqslant 1}\{v_{1,k}^{(i)}q_{1,k}(t)+v_{2,k}^{(i)}q_{2,k}(t)\},
\quad i=1,2
\end{equation}
and prove that the series converges. 
The moment problem  \eqref{moment} can be written as
\begin{equation}\label{moment Lambda1}
A_{1,k} V_{1,k} + A_{2,k}  V_{2,k} = F_k, ~~\mbox{ for all }k\geqslant 1,
\end{equation}
with for all $k\in\mathbb{N}^*$ 
\begin{equation}\label{Ak}
A_{1,k} = \left(
\begin{array}{cc}
 \widetilde f_k^{(1)} &  \widetilde f_k^{(2)}\\ 
f_k^{(1)} & f_k^{(2)} 
\end{array}
\right),~~~~
A_{2,k} =\left(
\begin{array}{cc}
-I_k(p,q) f_k^{(1)} &- I_k(p,q) f_k^{(2)}\\
0 & 0 
\end{array}
\right),
\end{equation}
 \begin{equation}\label{V_k}
V_{1,k}:=\left(
\begin{array}{cc}
v_{1,k}^{(1)} \\
v_{1,k}^{(2)}
\end{array}
\right), \quad 
V_{2,k}:= \left(
\begin{array}{cc} 
v_{2,k}^{(1)} \\
v_{2,k}^{(2)}
\end{array}\right)  
\end{equation}
and
\begin{equation}\label{Fk}
F_k = \left(
\begin{array}{c}
- e^{- k^2 T} \left(   y_{1,k}^0 - T I_k(p,q) y_{2,k}^0 \right)\\
- e^{- k^2 T}y_{2,k}^0 
\end{array}
\right).
\end{equation}

\subsection{Solving the moment problem}\label{sec: resol}

In this section, we will prove the 
null controllability of System \eqref{system primmal int} at any time $T$ 
when the supports of $p$ or $q$ intersects the control domain $\omega$ (Theorem 1.1). 
In \cite{gonzalez2010controllability}, the authors obtain 
the null controllability of System \eqref{system primmal int} at any time  
under Condition  \eqref{cond casc}, 
so we will not consider this case and  
we will always suppose that $\mathrm{Supp}~p\cap\omega\neq\varnothing$. 
This implies that there exists $x_0\in\omega$ such that $p(x_0)\neq0$. 
By continuity of $p$, we deduce that
$|p|>C$ in $\widetilde{\omega}$ for a positive constant $C$ 
and an open subinterval $\widetilde{\omega}$ of $\omega$. 

 \begin{dfntn}\label{def:equiv}
 Let $p_1,~p_2\in W^1_{\infty}(0,\pi)$ and $q_1,~q_2\in L^{\infty}(0,\pi)$. Consider the systems given for $i\in\{1,2\}$ by
 \begin{equation}\label{system Si}
 \left\{\begin{array}{l}
 \mathrm{For~given ~}y^0\in L^2(0,\pi)^2,~v\in L^2(Q_T),\\  \noalign{\smallskip}
 \mathrm{Find}~y:=(y_1,y_2)\in 
 L^2(0,T;H^1_0(0,\pi)^2)\cap \mathcal{C}([0,T];L^2(0,\pi)^2)\mathrm{~such~that}:\\  \noalign{\smallskip}
\begin{array}{ll} \partial_ty_1-\partial_{xx}y_1=\mathds{1}_{\omega}v&\mathrm{in~} Q_T,\\  \noalign{\smallskip}
\partial_ty_2-\partial_{xx}y_2+p_i(x)\partial_xy_1+q_i(x)y_1=0
&\mathrm{in~} Q_T,\\  \noalign{\smallskip}
y(0,\cdot)=y(\pi,\cdot)=0&\mathrm{on}~(0,T),\\  \noalign{\smallskip}
y(\cdot,0)=y^0&\mathrm{in}~(0,\pi).\tag{$\mathcal{S}_{i}$}
        \end{array}\end{array}
\right.
\end{equation}
We say that System $(\mathcal{S}_1)$ is \textit{equivalent} to System $(\mathcal{S}_2)$ 
if System $(\mathcal{S}_1)$ is null controllable at time $T$ if and only if System $(\mathcal{S}_2)$ is null controllable at time $T$. 
 \end{dfntn}
 
 Let us present the main technique used all along this section. 
 Suppose that System \eqref{system primmal int} is null controllable at time $T$. Let $v$ a control such that the solution $y$ to System \eqref{system primmal int} 
 verifies $y(T)=0$ in $(0,\pi)$ and $\omega_0:=(\alpha,\beta)$ a subinterval of $\omega=(a,b)$. 
Consider a function $\theta\in W^2_{\infty}(0,\pi)$ satisfying 
 \begin{equation}\label{cond gene theta}
  \left\{\begin{array}{ll}
       \theta \equiv\kappa_1&\mathrm{~in~}(0,\alpha),\\  \noalign{\smallskip}
       \theta\equiv \kappa_2&\mathrm{~in~} (\beta,\pi),\\  \noalign{\smallskip}
       \theta>\kappa_3&\mathrm{~in~}(0,\pi),
         \end{array}
\right.
 \end{equation}
with $\kappa_1,~\kappa_2,~\kappa_3\in\mathbb{R}^*_+$. Thus if we consider the change of unknown 
\begin{equation}\label{chgt var}
 \widehat{y}:=(\widehat{y}_1,y_2)
 \mathrm{~~with~~}
 \widehat{y}_1:=\theta^{-1}y_1,
\end{equation}
 then $\widehat{y}$ is solution in $L^2(0,T;H^1_0(0,\pi)^2)\cap \mathcal{C}([0,T];L^2(0,\pi)^2)$ to the system 
 \begin{equation}\label{system primmal int ch 1}
 \left\{\begin{array}{ll}
\partial_t\widehat{y}_1-\partial_{xx}\widehat{y}_1=\mathds{1}_{\omega}\widehat{v}&\mathrm{in~} Q_T,\\  \noalign{\smallskip}
\partial_ty_2-\partial_{xx} y_2+\widehat{p} \partial_x\widehat{y}_1+\widehat{q}\widehat{y}_1=0&\mathrm{in~} Q_T,\\   \noalign{\smallskip}
\widehat{y}(0,\cdot)=\widehat{y}(\pi,\cdot)=0&\mathrm{on}~(0,T),\\  \noalign{\smallskip}
\widehat{y}(\cdot,0)=\widehat{y}^0&\mathrm{in}~(0,\pi),
        \end{array}
\right.
\end{equation}
 where the initial condition is $\widehat{y}^0:= (\theta^{-1}y_1^0,y_2^0)\in L^2(0,\pi)^2$, 
 the control is $\widehat{v}:=-\partial_{xx}(\theta^{-1})y_1-2\partial_x(\theta^{-1})\partial_xy_1+\theta^{-1}\mathds{1}_{\omega}v\in L^2(Q_T)$ 
 and the coupling terms are given by $\widehat{p}:=p\theta$ and  $\widehat{q}:=p\partial_x\theta+q\theta$. 
Since $\theta$ is constant in $(0,\pi)\backslash\omega_0$, 
 we have $\mathrm{Supp}~\widehat{v}\subseteq \omega\times(0,T)$.
 Since $y$ is controlled, then $\widehat{y}$ also. 
 The converse is clearly true: starting from the controlled System \eqref{system primmal int ch 1} 
 the same process leads to the construction of a controlled solution 
 of System \eqref{system primmal int}. 
Thus through the change of unknown \eqref{chgt var}, following Definition \ref{def:equiv}, 
  Systems \eqref{system primmal int} and \eqref{system primmal int ch 1} are equivalent. 
 
The next main result of this section is Proposition \ref{prop intersec2} 
that will be introduced after some lemmas. 
The first of them is the following.

 \begin{lmm}\label{lem equiv q=0}Let $p\in W^1_{\infty}(0,\pi)\cap W^2_{\infty}(\omega)$ and $q\in L^{\infty}(0,\pi)\cap W^1_{\infty}(\omega)$ with $|p|>C$ in
 an open subinterval $\widetilde{\omega}$ of $\omega$ for a positive constant $C$. 
 There exists a subinterval $\omega_1:=(\alpha,\beta)\subset\widetilde{\omega}$ and a function $\theta\in W^2_{\infty}(0,\pi)$ satisfying \eqref{cond gene theta} such that 
 System \eqref{system primmal int} is equivalent to System \eqref{system primmal int ch 1} with 
$\widehat{q}\equiv 0$ in  $\omega_1$. 
Moreover, for all $\epsilon>0$, the interval  $\omega_1$ can be chosen in order to obtain
for all $k\in\mathbb{N}^*$
\begin{equation}\label{estim I Ik theta}
|I_k(p,q)-I_k(\widehat{p},\widehat{q})|\leqslant \varepsilon.
\end{equation} 
Consequently, taking the limit, we deduce that $|I(p,q)-I(\widehat{p},\widehat{q})|\leqslant \varepsilon$.
 \end{lmm}
\begin{proof}
Let  $\omega_1:=(\alpha,\beta)$ be an  interval strictly included in $\widetilde{\omega}:=(\widetilde{a},\widetilde{b})$ and $\theta\in W^2_{\infty}(0,\pi)$ satisfying
\begin{equation}\label{theta q=0}
\left\{\begin{array}{ll}
 p\partial_x\theta+q\theta=0&\mathrm{~in~}\omega_1,\\  \noalign{\smallskip}
 \theta\equiv 1&\mathrm{~in~}(0,\pi)\backslash\widetilde{\omega},\\  \noalign{\smallskip}
 |\theta|>C&\mathrm{~in~}(0,\pi),
\end{array}\right.\end{equation}
for a positive constant $C$.  
In the intervals $(\widetilde{a},\alpha]$ and $[\beta,\widetilde{b})$, we can take $\theta$ of class $\mathcal{C}^{\infty}$ in order to have $\theta\in W^2_{\infty}(0,\pi)$. 
Thus the function $\theta$ verifies \eqref{cond gene theta} 
and,  following the change of unknown described in \eqref{chgt var}, System \eqref{system primmal int} is equivalent to System \eqref{system primmal int ch 1} with $\widehat{q}\equiv 0$ in  $\omega_1$. 
The estimates in \eqref{estim I Ik theta} are obtained taking the interval $\omega_1$ small enough, so $\theta$  will be close to $1$.
\end{proof}

Let us first study System \eqref{system primmal int} in a particular case.
\begin{lmm}\label{prop intersec1}Consider $p\in W^1_{\infty}(0,\pi)\cap W^2_{\infty}(\omega)$ and $q\in L^{\infty}(0,\pi)\cap W^1_{\infty}(\omega)$.
 Let us suppose that $p\equiv C\in \mathbb{R}^*$ and $q\equiv0$  in an open subinterval $\widetilde{\omega}$ of $\omega$. 
Then   System \eqref{system primmal int} is equivalent to a system of the form \eqref{system primmal int ch 1} with coupling terms $\widehat{p},~\widehat{q}$ satisfying
 \begin{equation*}\begin{array}{c}
  |I_k(\widehat{p},\widehat{q})|>C/k^{6},~\forall k\in\mathbb{N}^*.
 \end{array}\end{equation*}
\end{lmm}

To prove this result we will need this lemma:

 \begin{lmm}\label{lemme suite}
Let  $(u_k)_{k\in\mathbb{N}^*}$ be a real sequence. 
Then there exists $\kappa\in\mathbb{R}^*_+$ such that for all $k\in\mathbb{N}^*$
\begin{equation*}
 |u_k+\kappa|\geqslant 1/k^2.
\end{equation*}
 \end{lmm}

  \begin{proof}[Proof of Lemma \ref{lemme suite}]
By contradiction let us suppose that for all $\kappa\in\mathbb{R}^*_+$ there exists $k\in\mathbb{N}^*$ such that 
 $|u_k+\kappa|< 1/k^2$.
Then 
\begin{equation}\label{preuve lemme 3.2}
 \mathbb{R}^*_+\subseteq\bigcup\limits_{k\in\mathbb{N}^*}(-u_k-1/k^2,-u_k+1/k^2).
\end{equation}
The convergence of 
 the series $\sum_{k\in \mathbb{N}^*}1/k^2$ implies that the measure of the set in the right-hand side in \eqref{preuve lemme 3.2} is finite and leads to the conclusion.
\end{proof}

\begin{proof}[Proof of Lemma \ref{prop intersec1}]
Let  $(\alpha,\beta)$ an open subinterval of $\widetilde{\omega}$ with $\alpha$ and $\beta$ to be determined later, 
$\kappa\in \mathbb{R}^*_+$ and  $\theta\in W^2_{\infty}(0,\pi)$ satisfying
 \begin{equation}\label{def theta p=C}
  \left\{\begin{array}{ll}
       \theta \equiv 1&\mathrm{~in~}(0,\pi)\backslash (\alpha,\beta),\\ \noalign{\smallskip}
       \theta\equiv 1+ \kappa\xi&\mathrm{~in~} (\alpha,\beta),       
         \end{array}
\right.
 \end{equation}
 where
 \begin{equation}\label{def xi}
 \xi(x):=\left[\sin\left(\frac{\pi(x-\alpha)}{\beta-\alpha}\right)\right]^2
 \mathrm{~for~all~}x\in(\alpha,\beta).
 \end{equation}
 In particular, we have  $\theta\geqslant1$ in $(0,\pi)$.
Let $k\in\mathbb{N}^*$, $\widehat{y}_1:=\theta^{-1}y_1$ and  $\widehat{y}:=(\widehat{y}_1,y_2)$ 
the solution  
to System \eqref{system primmal int ch 1}. 
For System \eqref{system primmal int ch 1} the quantity $I_k$ defined in the introduction is given by 
 \begin{equation*}\begin{array}{rcl}
  I_k(\widehat{p},\widehat{q})&=&\displaystyle\int_0^{\pi}
  \{\widehat{q}-\frac{1}{2}\partial_x\widehat{p}\}\varphi_k^2\vspace*{2mm}\\
  &=&I_k(p,q)+\kappa J_k,
 \end{array}\end{equation*}
 with $\widehat{p},~\widehat{q}$ given by $\widehat{p}:=p\theta$ and  $\widehat{q}:=p\partial_x\theta+q\theta$ and $J_k$ defined by
 \begin{equation*}
   J_k:=\frac{1}{2}\displaystyle\int_{\alpha}^{\beta}\partial_x(\xi)\varphi_k^2.
 \end{equation*}
Then, after a simple calculation, we obtain
  \begin{equation}\label{calcul Jk 1}\begin{array}{rcl}
 J_k 
&=&\dfrac{\frac{2\pi}{(\beta-\alpha)^2}}{(2k+\frac{2\pi}{\beta-\alpha})(2k-\frac{2\pi}{\beta-\alpha})}\sin(k(\beta+\alpha))\sin(k(\beta-\alpha)). 
 \end{array}\end{equation}
 Let $n\in\mathbb{N}^*$ large enough such that 
  $\frac{a}{n}<\frac{b}{n+1}$.
There exists $\ell$ an algebraic number of order two, \textit{i.e.} 
a root of a  polynomial of degree 2 with integer coefficients, 
 satisfying 
 \begin{equation*}
  \frac{a}{n}<\ell<\frac{b}{n+1}
 \hspace*{5mm} \mathrm{and}\hspace*{5mm} \ell\neq \frac{\pi}{j}\hspace*{5mm} 
 \mathrm{for~all~}j\in\mathbb{N}^*,
 \end{equation*}
 since the set of such numbers is dense in $\mathbb{R}$. 
Let us take $\alpha:=n\ell$ and $\beta:=(n+1)\ell$. Thus $\alpha,~\beta\in(a,b)$,
\begin{equation}\label{calcul Jk 2}
k(\beta+\alpha)=k(2n+1)\ell\hspace*{5mm} \mathrm{and}\hspace*{5mm} 
k(\beta-\alpha)=k\ell.
\end{equation}
Moreover
\begin{equation*}
 \left|2k+\frac{2\pi}{\beta-\alpha}\right|\times
\left|2k-\frac{2\pi}{\beta-\alpha}\right|<Rk^2,
\end{equation*}
with $R>0$. Since $\ell$ is an algebraic number of order two, using diophantine approximations it can be proved that
\begin{equation}\label{calcul Jk 3}
 \inf\limits_{j\geqslant 1}(j|\sin(j\ell)|)\geqslant \gamma,
\end{equation}
for a positive constant $\gamma$ (see \cite[Ine. (5.13)]{Ammar-Khodja2015}).
The expressions \eqref{calcul Jk 1}-\eqref{calcul Jk 3} give
\begin{equation}\label{estim J_k}
 |J_k|\geqslant \dfrac{2\pi}{(\beta-\alpha)^2}\frac{\gamma^2}{R(2n+1)k^4}
\end{equation}
for all $k\in \mathbb{N}^*$. 
Using Lemma \ref{lemme suite}, there exists $\kappa\in\mathbb{R}^*_+$ satisfying
\begin{equation*}
 \left|\frac{I_k(p,q)}{J_k}+\kappa\right|\geqslant 1/k^2.
\end{equation*}
 Combining the last inequality with Estimate \eqref{estim J_k},
 \begin{equation*}
 |I_k(\widehat{p},\widehat{q})|=\left|I_k(p,q)+\kappa J_k\right|\geqslant |J_k|/k^2\geqslant C/k^{6}.
\end{equation*}
\end{proof}

The next lemma is proved in \cite[Lem. 5.1]{Ammar-Khodja2015}.

\begin{lmm}\label{lemme fk}
 There exist functions $f^{(1)}$, $f^{(2)}\in L^2(0,\pi)$ satisfying  $\mathrm{Supp}~f^{(1)},~\mathrm{Supp}~f^{(2)}\subseteq\omega$ 
 and such that for all $k\in\mathbb{N}^*$
 \begin{equation}\label{estim fik hat fik}
\mathrm{min}\{|f_k^{(1)}|,|f_k^{(2)}|\}\geqslant \dfrac{C}{k^3}
\hspace*{5mm}\mbox{and}\hspace*{5mm}
|B_k|:=|\widehat{f}_k^{(1)}{f}_{k}^{(2)}-\widehat{f}_k^{(2)}{f}_{k}^{(1)}|\geqslant \dfrac{C}{k^{5}},
 \end{equation}
where for $i\in\{1,2\}$ the terms $f_{k}^{(i)}$ and $\widehat{f}_{k}^{(i)}$ are given by
\begin{equation}\label{def hat fik}
 f_{k}^{(i)}:=\int_0^{\pi}f^{(i)}(x)\varphi_k(x)dx\hspace*{5mm}\mathrm{and}\hspace*{5mm}
 \widehat{f}_{k}^{(i)}:=\int_0^{\pi}f^{(i)}(x)\cos(kx)dx.
\end{equation}
\end{lmm}

With the help of Lemma \ref{lemme fk}, we deduce the following proposition:

\begin{prpstn}\label{prop intersec2}Consider $p\in W^1_{\infty}(0,\pi)\cap W^2_{\infty}(\omega)$ and $q\in L^{\infty}(0,\pi)\cap W^1_{\infty}(\omega)$.
 Let us suppose that $|p|>C$ in   an open subinterval  $\widetilde{\omega}$ of $\omega$ for a positive constant $C$. 
 Then System \eqref{system primmal int} is equivalent to a system of the form \eqref{system primmal int ch 1} 
 with coupling terms $\widehat{p},~\widehat{q}$ satisfying Condition \eqref{cond approx2}, $T_0(\widehat{p},\widehat{q})=0$ and 
\begin{equation}\label{cond prop intersec2}
\left| \det A_{1,k}  \right| \ge \frac{C_1}{k^{7}} \left| I_{a,k} (\widehat{p},\widehat{q}) \right| 
- \frac{C_2}{k} \left| I_k (\widehat{p},\widehat{q}) \right|  ~~\forall k\in\mathbb{N}^*,
\end{equation} 
 where $C_1$ and $C_2$ are two positive constants independent on $k$ 
 (the notion of equivalent systems is defined at the beginning of Section 3.2). 
\end{prpstn}

\begin{proof}
Using Lemma \ref{lem equiv q=0}, 
without loss of generality, we can suppose that $q\equiv0$ and $|p|>C$ in a subinterval $\widehat{\omega}$ of $\widetilde{\omega}$ for a positive constant $C$. 
 If $\partial_xp\equiv 0$ in $\widehat{\omega}$, Lemma \ref{prop intersec1} leads to 
 \begin{equation*}
 | I_k(p,q)|\geqslant C/k^6,~~\forall k\in\mathbb{N}^*,
 \end{equation*}
which implies that Condition \eqref{cond approx2} is satisfied and the right-hand side of inequality \eqref{cond prop intersec2} is negative  
for some appropriate constants $C_1$ and $C_2$. 
 Otherwise, let $(\alpha,\beta)\subseteq\widehat{\omega}$ such that $\partial_xp>C$   in $(\alpha,\beta)$ 
 or $\partial_xp<-C$ in $(\alpha,\beta)$ for a positive constant $C$.  
The rest of the  proof is divided into three steps: 

\textbf{Step 1:} 
If $I(p,q):=\int^{\pi}_0\{q-\frac12\partial_xp\}=0$, 
we will prove in this step 
that System \eqref{system primmal int} is equivalent to a system with coupling terms 
$\widehat{p},~\widehat{q}$ satisfying $  I(\widehat{p},\widehat{q})\neq0$. 
Assume that $I(p,q)=0$ and  consider $\theta\in W^2_{\infty}(0,\pi)$ 
defined in \eqref{def theta p=C}, with $\kappa:=1$. We remark that $|\theta|\geqslant 1$. 
If we consider the change of unknown described in \eqref{chgt var}, 
then for all $k\in\mathbb{N}^ *$, 
using the definition of $I_k$, we obtain
 \begin{equation*}\begin{array}{rcl}
  I_k(\widehat{p},\widehat{q})&=&
  I_k(p,q)+\displaystyle\int_{\alpha}^{\beta}\{\frac12\partial_x(\xi)p
  -\frac12\xi\partial_x(p)\}\varphi_k^2dx\\  \noalign{\smallskip}
   & =&I_k(p,q)+ J_{k}(p,q),
\end{array} \end{equation*}
where
\begin{equation*}\begin{array}{rcl}
 J_{k}(p,q)&=&\dfrac{1}{2\pi}\displaystyle\int_{\alpha}^{\beta} \{\partial_x(\xi)p-\xi\partial_x(p)\}\{1-\cos(2kx)\}dx
\\ \noalign{\smallskip}
 &\substack{\longrightarrow\\k\rightarrow \infty}&
 -\dfrac{1}{\pi}\displaystyle\int_{\alpha}^{\beta}\xi\partial_x(p)dx=:J(p,q).
\end{array}\end{equation*}
Using the definition of $\xi$ given in \eqref{def xi}, we get
\begin{equation*}\begin{array}{rcl}
 |J(p,q)|&\geqslant&\dfrac{1}{\pi}
 \inf\limits_{(\alpha,\beta)}|\partial_xp|
 \displaystyle\int_{\alpha}^{\beta}\sin^2\left(\frac{\pi(x-\alpha)}{\beta-\alpha}\right)dx\\ \noalign{\smallskip}
 &=&\dfrac{1}{2\pi}
 \inf\limits_{(\alpha,\beta)}|\partial_xp|
 \displaystyle\int_{\alpha}^{\beta}\{1-\cos\left(\frac{2\pi(x-\alpha)}{\beta-\alpha}\right)\}dx\\ \noalign{\smallskip}
&=& \dfrac{(\beta-\alpha)}{2\pi}
 \inf\limits_{(\alpha,\beta)}|\partial_xp|
 \neq0.
 \end{array}\end{equation*}
We  recall that $I_k(p,q)\rightarrow I(p,q)=0$. Thus, we obtain 
$I_k(\widehat{p},\widehat{q})\rightarrow I(\widehat{p},\widehat{q})=J(p,q)\neq0$.

\textbf{Step 2:}
We will show in this second step that 
System \eqref{system primmal int} is equivalent to a system with coupling terms $\widehat{p},~\widehat{q}$ such that  $ |I_k(p,q)|>C>0\mbox{ for all }k\in\mathbb{N}^*\mbox{ satisfying }p\varphi_k\mbox{ non-constant}$.  
In view of Step 1, we can assume that $I(p,q)\neq0$. 
Using Lemma \ref{lem equiv q=0}, up to the change of unknown \eqref{chgt var}  we can also suppose that $q\equiv 0$ in  an open subinterval $\widehat{\omega}$ of $\widetilde{\omega}$. 
Moreover, by \eqref{estim I Ik theta}, the function $\theta$ and $\widehat{\omega}$ can be chosen  in order to keep  the quantity $I$ different of zero. 
Let $(\alpha,\beta)\subseteq\widehat{\omega}$ such that $|p|>C>0$ in $(\alpha,\beta)$.
 Since $I(p,q)\neq0$ and $I_k(p,q)\rightarrow I(p,q)$, there exists $k_0\in\mathbb{N}^*$ 
such that  $|I_k(p,q)|>C$ for a  constant $C>0$ and all $k\geqslant k_0$.
Let us define the set 
$$S_0:=\{k\in\mathbb{N}^*:I_k(p,q)=0
\hspace*{5mm} \mathrm{and}\hspace*{5mm} 
p\varphi_k\mathrm{~non-constant~}\mathrm{~in~}(\alpha,\beta)\}$$
and $M:=\#S_0<\infty$. 
 Let $\theta\in W^2_{\infty}(0,\pi)$ satisfying
 \begin{equation*}
  \left\{\begin{array}{ll}
       \theta=1+\sum_{m=1}^{M}\xi_m,~ |\theta|>C>0,&\\  \noalign{\smallskip}
       \xi_m\in W^2_{\infty}(0,\pi),~
       \mathrm{Supp}~\xi_m\subseteq(\alpha,\beta),&\mathrm{for~all~} m\in\{1,...,M\},
         \end{array}
\right.
 \end{equation*}
where $\xi_1,...,\xi_M$ are to be determined. 
Again, if we consider the change of unknown \eqref{chgt var}, 
then for all $k\in\mathbb{N}^ *$, 
using the definition of $I_k$, we obtain
 \begin{equation*}\begin{array}{rcl}
  I_k(\widehat{p},\widehat{q})&=&
  I_k(p,q)
+\sum\limits_{m=1}^{M}\displaystyle\int_{\alpha}^{\beta}\{\frac12\partial_x(\xi_m)p-\frac12\xi_m\partial_x(p)\}\varphi_k^2dx\\  \noalign{\smallskip}
   & =:&I_k(p,q)+\sum\limits_{m=1}^MJ_{m,k}(p,q).
\end{array} \end{equation*}
The goal is to choose the functions $\xi_1,...,\xi_M$ such that 
for a constant $C>0$ we have $|I_k(\widehat{p},\widehat{q})|>C$ 
for all $k\in\mathbb{N}^*$ satisfying $p\varphi_k$ non-constant in $(\alpha,\beta)$. 
 We will construct $\xi_1,...,\xi_{M}$ from $\xi_1$ until $\xi_M$. 

Let $k\in S_{0}$  
 and  consider  $(f_1,\xi_1)\in  W^1_{\infty}(\alpha,\beta)\times W^2_{\infty}(\alpha,\beta)$ 
 a solution to
 \begin{equation*}
  \left\{\begin{array}{ll}
       \dfrac12\partial_x(\xi_1)p-\dfrac12\xi_1\partial_x(p)=f_1&\mathrm{in}~(\alpha,\beta),\\ \noalign{\smallskip}
       \xi_1(\alpha)=\xi_1(\beta)
       =\partial_x\xi_1(\alpha)=\partial_x\xi_1(\beta)=0.&
         \end{array}
\right.
 \end{equation*}
This system is equivalent to
 \begin{equation*}
  \left\{\begin{array}{l}
       \xi_1(x)=p(x)\displaystyle\int_{\alpha}^x\frac{2f_1(s)}{p^2(s)}ds,\mathrm{~~~~for~all~}x\in(\alpha,\beta),\\ \noalign{\smallskip}
       \displaystyle\int_{\alpha}^{\beta}\frac{2f_1(s)}{p^2(s)}ds=0, f_1(\alpha)=f_1(\beta)=0.
                \end{array}
\right.
 \end{equation*}
We remark that we need that $p\in W^2_{\infty}(\alpha,\beta)$. Finding a function $f_1$ satisfying
  \begin{equation}\label{constr f J0}
         f_1(\alpha)=f_1(\beta)=0,
        \displaystyle \int_{\alpha}^{\beta}\frac{2f_1(s)}{p^2(s)}ds=0 
        \hspace*{5mm}\mathrm{and}\hspace*{5mm}
       J_{1,k}(p,q)=  \displaystyle\int_{\alpha}^{\beta}f_1(s)\varphi_k^2(s)ds\neq0,
 \end{equation}
 is equivalent 
 to finding a function $g:=2f_1/p^2$ satisfying 
     \begin{equation*}
         g_1(\alpha)=g_1(\beta)=0,
       \displaystyle \int_{\alpha}^{\beta}g_1(s)ds=0
       \hspace*{5mm}\mathrm{and}\hspace*{5mm}
       \displaystyle \int_{\alpha}^{\beta}g_1(s)p^2(s)\varphi_k^2(s)ds\neq0.
 \end{equation*}     
 Let $\kappa_1\in\mathbb{R}$ and  define for all $j\in\mathbb{N}^*$ and all $x\in(\alpha,\beta)$
\begin{equation*}
 g_{1,j}(x):=\kappa_1\sin\left(\frac{2\pi j(x-\alpha)}{\beta-\alpha}\right).
\end{equation*}
Using the fact that $p\varphi_k$ is non-constant in $(\alpha,\beta)$, 
without loss of generality, we can suppose that 
$$\varphi_k\left(\alpha+\frac{\beta-\alpha}{4}\right)p\left(\alpha+\frac{\beta-\alpha}{4}\right)
\neq\varphi_k\left(\alpha+\frac{3(\beta-\alpha)}{4}\right)p\left(\alpha+\frac{3(\beta-\alpha)}{4}\right),$$
otherwise we adapt the interval $(\alpha,\beta)$ at the beginning of Step 2. 
We  deduce that the function $h_k$ of $L^2(\alpha,\alpha+(\beta-\alpha)/2)$ 
defined by
\begin{equation*}\begin{array}{cccl}
h_k(s):=p^2(s)\varphi_k^2(s)-p^2(\beta+\alpha-s)\varphi_k^2(\beta+\alpha-s)
\end{array}
 \end{equation*}
is not equal to zero in $(\alpha,\alpha+(\beta-\alpha)/2)$. 
Since $(g_{1,j})_{j\in\mathbb{N}^*}$ is a Riesz basis of  $L^2(\alpha,\alpha+(\beta-\alpha)/2)$, 
there exists $j_1\in \mathbb{N}^*$ such that
\begin{equation*}
\displaystyle\int_{\alpha}^{\alpha+(\beta-\alpha)/2}g_{1,j_1}(s)
\left[p^2(s)\varphi_k^2(s)-p^2(\beta+\alpha-s)\varphi_k^2(\beta+\alpha-s)
\right]ds
\neq 0.
 \end{equation*}
Moreover, using the fact that 
 $g_{1,j_1}(s)=g_{1,j_1}(\beta+\alpha-s)~\forall s\in (\alpha,\alpha+(\beta-\alpha)/2)$, 
we have
\begin{equation*}
\displaystyle\int_{\alpha}^{\alpha+(\beta-\alpha)/2}g_{1,j_1}(s)p^2(s)\varphi_k^2(s)ds
\neq -\displaystyle\int_{\alpha+(\beta-\alpha)/2}^{\beta}
g_{1,j_1}(s)p^2(s)\varphi_k^2(s)ds.
 \end{equation*}
Thus 
 \begin{equation*}
 \displaystyle\int_{\alpha}^{\beta}g_{1,j_1}(s)p^2(s)\varphi_k^2(s)ds
\neq 0.
 \end{equation*}
 Plugging  $g_1:=g_{1,j_1}$ and $f_1:=\dfrac{g_{1,j_1}p^2}{2}$ in \eqref{constr f J0}, 
 we obtain 
 \begin{equation*}
  J_{1,k}(p,q)=\frac{\kappa_1}{2}\displaystyle\int_{\alpha}^{\beta}\sin\left(
  \frac{2\pi j_1(s-\alpha)}{\beta-\alpha}\right)p(s)^2\varphi_k(s)^2ds\neq0.
 \end{equation*}
 We have also for all $j\in\mathbb{N}^*$
  \begin{equation*}
  J_{1,j}(p,q)=\frac{\kappa_1}{2}\displaystyle\int_{\alpha}^{\beta}\sin\left(
  \frac{2\pi j_1(s-\alpha)}{\beta-\alpha}\right)p(s)^2\varphi_j(s)^2ds.
 \end{equation*}
We fix $\kappa_1$ in order to have 
$$\sup\limits_{i\in\mathbb{N}^*}|J_{1,i}(p,q)|
\leqslant \frac12\inf\limits_{i\in\mathbb{N}^*\backslash S_0}\left|I_i(p,q)\right|.$$

Let $m\in\{2,...,M\}$ and let us assume that $\xi_1,...,\xi_{m-1}$ are already constructed. 
 Consider  the set 
 $$S_{m-1}:=\{k\in\mathbb{N}^*:I_k(p,q)+\sum_{j=1}^{m-1}J_{j,k}(p,q)=0
 \mathrm{~and~}p\varphi_k\mathrm{~non-constant~in~}(\alpha,\beta)\}.$$  
 If $S_{m-1}=\varnothing$,  
then we take $\xi_m=0$ in $(0,\pi)$. 
 Otherwise, let $k\in S_{m-1}$ 
 and  consider  $(f_m,\xi_m)\in  W^1_{\infty}(\alpha,\beta)\times W^2_{\infty}(\alpha,\beta)$ 
 a solution to
 \begin{equation*}
  \left\{\begin{array}{ll}
       \dfrac12\partial_x(\xi_m)p-\dfrac12\xi_m\partial_x(p)=f_m&\mathrm{in}~(\alpha,\beta),\\ \noalign{\smallskip}
       \xi_m(\alpha)=\xi_m(\beta)
       =\partial_x\xi_m(\alpha)=\partial_x\xi_m(\beta)=0.&
         \end{array}
\right.
 \end{equation*}
This system is equivalent to
 \begin{equation*}
  \left\{\begin{array}{l}
       \xi_m(x)=p(x)\displaystyle\int_{\alpha}^x\frac{2f_m(s)}{p^2(s)}ds,\mathrm{~~~~for~all~}x\in(\alpha,\beta),\\ \noalign{\smallskip}
       \displaystyle\int_{\alpha}^{\beta}\frac{2f_m(s)}{p^2(s)}ds=0, f_m(\alpha)=f_m(\beta)=0.
         \end{array}
\right.
 \end{equation*}
 Let $\kappa_m>0$. Again, there exists $j_m\in\mathbb{N}^*$ such that the function $f_m$ given 
 for all $x\in (\alpha,\beta)$ by
 \begin{equation*}
  f_m(x):=\frac{\kappa_m}{2}\sin\left(
  \frac{2\pi j_m(x-\alpha)}{\beta-\alpha}\right)p(x)^2
 \end{equation*}
 is solution to this system. 
 Then, we obtain
  \begin{equation*}
  J_{m,j}(p,q)=\frac{\kappa_m}{2}\displaystyle\int_{\alpha}^{\beta}\sin\left(
  \frac{2\pi j_m(s-\alpha)}{\beta-\alpha}\right)p(s)^2\varphi_j(s)^2ds.
 \end{equation*} 
The last quantity is different of zero for $j=k$. 
Let us fix $\kappa_m$ in order to have 
$$\sup\limits_{i\in\mathbb{N}^*}|J_{m,i}(p,q)|
\leqslant \frac12\inf\limits_{i\in\mathbb{N}^*\backslash S_{m-1}}
\left|I_i(p,q)+\sum_{j=1}^{m-1}J_{j,i}(p,q)\right|.$$
Thus, after constructing the functions $\xi_1,...,\xi_M$, the obtained functions $\widehat{p}$ and $\widehat{q}$ are such that 
$$|I_k(\widehat{p},\widehat{q})|>C 
\mbox{ for all }k\in\mathbb{N}^*\mbox{ satisfying }\widehat{p}\varphi_k
\mbox{ non-constant in }(\alpha,\beta),$$
where $C$ is a positive constant which does not depend on $k$.

\textbf{Step 3:} 
Finally, in this third step, we will prove  that System \eqref{system primmal int} 
 is equivalent to a system 
 satisfying $T_0(\widehat{p},\widehat{q})=0$ and Conditions \eqref{cond approx2} 
 and \eqref{cond prop intersec2}.
In view of Step 2, we can assume that 
$$|I_k(p,q)|>C 
\mbox{ for all }k\in\mathbb{N}^*\mbox{ satisfying }p\varphi_k
\mbox{ non-constant in }(\alpha,\beta),$$
where $C$ is a positive constant which does not depend on $k$. 
If $|I_k(p,q)|>C_0$ for all $k\in\mathbb{N}^*$ and a constant $C_0>0$, 
 then Condition \eqref{cond approx2} is satisfied and the right-hand side of inequality \eqref{cond prop intersec2} is negative  
for some appropriate constants $C_1$ and $C_2$. 
Let us now suppose that, for a $m\in\mathbb{N}^*$, we have
 \begin{equation}\label{sous cond inter}
 I_m(p,q)=0 \mbox{ and }p\varphi_m\mbox{ constant in }(\alpha,\beta).
 \end{equation}
Again, using Lemma \ref{lem equiv q=0}, up to the change of unknown \eqref{chgt var} described at the beginning of the section we can also suppose that $q\equiv 0$ in  a subinterval $(\alpha,\beta)$ of $\widetilde{\omega}$. 
Moreover, using \eqref{estim I Ik theta},
this change of unknown can be chosen in order to keep the property: 
$|I_k(p,q)|>C>0\mbox{ for all }k\in \mathbb{N}^*\backslash \{m\}$. 
Let $m\in\mathbb{N}^*$ such that $I_m(p,q)=0$ and  $p\varphi_m$ is constant in $(\alpha,\beta)$, 
otherwise we argue as in Step 2. 
Let  $\theta\in W^2_{\infty}(0,\pi)$ satisfying
 \begin{equation*}
  \left\{\begin{array}{l}
       \theta=1+\xi~~~~~\mathrm{in}~(0,\pi),\\  \noalign{\smallskip}
       \xi\in W^2_{\infty}(0,\pi), |\theta|>C>0,\\  \noalign{\smallskip}
       \xi\equiv\xi_{\alpha}\in\mathbb{R}^*_+\mathrm{~in~}(0,\alpha),\\  \noalign{\smallskip}
       \xi\equiv0\mathrm{~in~}(\beta,\pi),\\  \noalign{\smallskip}
         \end{array}
\right.
 \end{equation*}
 Again, if we consider the change of unknown described in \eqref{chgt var}, 
 then for all $k\in \mathbb{N}^*$
 \begin{equation*}\begin{array}{rcl}
  I_k(\widehat{p},\widehat{q})&=&
  I_k(p,q)+\displaystyle\int_{0}^{\beta}\{\frac12\partial_x(\xi)p+\xi q-\frac12\xi\partial_x(p)\}\varphi_k^2dx\\  \noalign{\smallskip}
   & =:&I_k(p,q)+J_{k}(p,q).
\end{array} \end{equation*}
We will distinguish the cases $I_{\alpha,m}(p,q)=0$ and $I_{\alpha,m}(p,q)\neq0$ (see \eqref{def Ik} for the definition of this quantity)
for the new control domain $\omega:=(\alpha,\beta)$.
\begin{enumerate}
 \item[Case 1:] Assume that $I_{\alpha,m}(p,q)=0$. 
  Let  $(\xi,h)\in W^2_{\infty}(\alpha,\beta)\times W^1_{\infty}(\alpha,\beta)$ 
  be a solution  to the system
 \begin{equation*}
  \left\{\begin{array}{ll}
       \frac12\partial_x(\xi)p-\frac12\xi\partial_x(p)=h&\mathrm{in}~(\alpha,\beta),\\  \noalign{\smallskip}
       \xi(\beta)=\partial_x\xi(\alpha)=\partial_x\xi(\beta)=0,~  \xi(\alpha)=\xi_{\alpha}\in\mathbb{R}^*.
         \end{array}
\right.
 \end{equation*}
This system is equivalent to
 \begin{equation*}
  \left\{\begin{array}{l}
       \xi(x)=-p(x)\int^{\beta}_x\frac{2h(s)}{p^2(s)}ds,\mathrm{~~~~~~for~all~}x\in (\alpha,\beta),\\  \noalign{\smallskip}
       \int_{\alpha}^{\beta}\frac{2h(s)}{p^2(s)}ds=\frac{-\xi_{\alpha}}{p(\alpha)},~
       h(\alpha)=\dfrac{-\xi_{\alpha}\partial_xp(\alpha)}{2}, h(\beta)=0.
         \end{array}
\right.
 \end{equation*}
Taking into account that $I_{\alpha,m}(p,q)=0$, $q\equiv0$ in $(\alpha,\beta)$ 
and $p\varphi_m\equiv\gamma$  in $(\alpha,\beta)$ for a $\gamma\in\mathbb{R}^*$, one gets
 \begin{equation*}
  J_{m}(p,q)=\xi_{\alpha}\displaystyle\int_0^{\alpha}(q-\frac12\partial_x(p))\varphi_m^2dx
  +\dfrac{\gamma^2}{2}\displaystyle\int_{\alpha}^{\beta}\partial_x\left(\frac{\xi}{p}\right)dx
  =-\frac{\gamma^2\xi_{\alpha}}{2p(\alpha)}\neq0.
 \end{equation*}
Let $\xi_{\alpha}$ and $h$ be such that
 $$\sup\limits_{k\in\mathbb{N}^*}|J_{k}(p,q)|
\leqslant \frac12\inf\limits_{k\in\mathbb{N}^*\backslash\{m\}}|I_k(p,q)|.$$
 Then $|I_k(\widehat{p},\widehat{q})|>C$ for all $k\in\mathbb{N}^*$ and a positive constant $C$. 
 Thus Condition \eqref{cond approx2} is satisfied and the right-hand side of inequality \eqref{cond prop intersec2} is negative  
for some appropriate constants $C_1$ and $C_2$. 
  \item[Case 2:] Let us now assume that $I_{\alpha,m}(p,q)\neq0$. Then Condition \eqref{cond approx2} is verified. 
 We recall that, in the moment problem described in the last section, we have 
 \begin{equation*}
  \mathrm{det}~A_{1,m}=\widetilde{f}_{m}^{(1)}{f}_{m}^{(2)}-\widetilde{f}_{m}^{(2)}{f}_{m}^{(1)},
 \end{equation*}
where $f_{m}^{(1)}$, $f_{m}^{(2)}$, $\widetilde{f}_{m}^{(1)}$ and $\widetilde{f}_{m}^{(2)}$ 
are given in \eqref{def fik}. 
Since $p\varphi_m$ is constant in $(\alpha,\beta)$, 
the function $\psi^*_m$ of Proposition \ref{prop base} reads for all $x\in(\alpha,\beta)$
\begin{equation*}
\left\{\begin{array}{l}
\begin{array}{rcl}
\psi_m^*(x)&=&\alpha_m^*\varphi_m-\dfrac{1}{m}\displaystyle\int_0^{\alpha}\sin(m(x-\xi))
 [\partial_x(p(\xi)\varphi_m(\xi))-q(\xi)\varphi_m(\xi)]d\xi\\ \noalign{\smallskip}
 &=&\tau_m\varphi_m(x)-\sqrt{\dfrac{\pi}{2}}\dfrac{1}{m}I_{\alpha,m}(p,q)\cos(mx),
\end{array},\\
\tau_m:=\alpha_m^*-\sqrt{\frac{\pi}{2}}\frac{1}{m}\displaystyle\int_0^{\alpha}\cos(m\xi)
[\partial_x(p(\xi)\varphi_m(\xi))-q(\xi)\varphi_m(\xi)]d\xi.
\end{array}\right.\end{equation*}
We deduce that  
\begin{equation*}
 \mathrm{det} A_{1,m}=-\sqrt{\frac{\pi}{2}}\frac{1}{m}I_{\alpha,m}(p,q)(\widehat{f}_{m}^{(1)}{f}_{m}^{(2)}-\widehat{f}_{m}^{(2)}{f}_{m}^{(1)}),
\end{equation*}
where $\widehat{f}_{m}^{(1)}$ and $\widehat{f}_{m}^{(2)}$ are given in Lemma \ref{lemme fk}. 
Using  Lemma \ref{lemme fk}, we obtain $\mathrm{det} A_{1,m}\neq0$. 
Thus, for $C_1$ small enough \eqref{cond prop intersec2} is true for $k=m$ 
and, for all $k\neq m$, the right-hand side of \eqref{cond prop intersec2} is negative for $C_2$ be enough.
\end{enumerate}
We conclude this proof remarking that, in each case, there exists $C>0$ and $k_0\in\mathbb{N}^*$ such that, 
for all $k\geqslant k_0$, we have 
$|I_k(\widehat{p},\widehat{q})|\geqslant C/k^6$,
which implies that $T_0(\widehat{p},\widehat{q})=0$.
\end{proof}

We recall that  $T_0(p,q)$ is given by \eqref{def T_0}. 
Before  proving Theorem \ref{theo null int}, we will establish the following proposition 
which is true even in the case where the coupling region and 
the control domain are disjoint. 

\begin{prpstn}\label{prop resol}
Assume that  Conditions \eqref{cond approx2} and \eqref{cond prop intersec2} hold and $T>T_0(p,q)$.  \\
Then System \eqref{system primmal int} is null controllable at time $T$.
\end{prpstn}

\begin{proof}
We will use the same  strategy than \cite{Ammar-Khodja2015}. 
Let $\varepsilon>0$. Using  the definition of the minimal time $T_0 (p,q)$ in \eqref{def T_0}, 
there exists a positive integer $k_\varepsilon$ for which 
\begin{equation}\label{kepsilon}
\min \left\{\log\left|I_{a,k} (p,q )\right|^{-1}, \log\left|I_{k} ( p,q ) \right|^{-1} \right\} 
< k^2(T_0 (p,q) + \varepsilon), \quad \forall k > k_\varepsilon.
\end{equation}
The goal is to solve the moment problem described in Section \ref{moment problem}. 
We recall that we look for a control $v$ of the form  \eqref{contv} and \eqref{contvi} 
with $ f^{(1)}$ and $ f^{(2)}$ defined in Lemma \ref{lemme fk}. 
We will solve the moment problem \eqref{moment Lambda1} 
depending on whether $k$ belongs to $\Lambda_{1}$, $\Lambda_{2}$ or $\Lambda_3$, where 
\begin{equation*}
 \left\{\begin{array}{lcl}
       \Lambda_{1}&:=&\{k\in \mathbb{N}^* : I_{k}(p,q)\neq0, ~I_{a,k}(p,q)\neq0\},\\\noalign{\smallskip}
       \Lambda_{2}&:=&\{k\in \mathbb{N}^* : I_k(p,q)\neq0, ~ I_{a,k}(p,q)=0\},\\\noalign{\smallskip}
        \Lambda_{3}&:=&\{k\in \mathbb{N}^* :  I_k(p,q)=0, ~ I_{a,k}(p,q)\neq0\}.
        \end{array}
\right.
\end{equation*}


\textbf{Case 1 :} Consider the case  $k \in \Lambda_{1}$ with $k \le k_{\varepsilon}$. \\
Let us take $v^{(2)}_{1,k}=v^{(2)}_{2,k}=0$. 
The moment problem \eqref{moment Lambda1} becomes
\begin{equation*}
\left\{\begin{array}{l}
\widetilde{f}_{k}^{(1)}v^{(1)}_{1,k}-I_k(p,q)f_k^{(1)}v^{(1)}_{2,k}=
- e^{- k^2 T} \left(  y_{1,k}^0 - T I_k(p,q) y_{2,k}^0\right), \\ \noalign{\smallskip}
f_k^{(1)}v_{1,k}^{(1)}=- e^{- k^2 T}y_{2,k}^0 .
\end{array}\right.
\end{equation*}
Since $I_k(p,q)\neq0$ and using the estimate  of $f_k^{(1)}$ and $f_k^{(2)}$ in Lemma 
\ref{lemme fk}, the last system has a unique solution 
\begin{equation}\label{sol moment case 1}
\left\{\begin{array}{l}
v_{1,k}^{(1)}=- e^{- k^2 T}\frac{y_{2,k}^0}{f_k^{(1)}},\\ \noalign{\smallskip}
v^{(1)}_{2,k}=
 \frac{e^{- k^2 T}}{I_k(p,q)f_k^{(1)}} 
\left( y_{1,k}^0-T I_k(p,q) y_{2,k}^0
-\widetilde{f}_{k}^{(1)}\frac{y_{2,k}^0}{f_k^{(1)}}\right).
\end{array}\right.
\end{equation}
Moreover, since the set of the $k$ considered in this case is finite, we get the inequality
\begin{equation}\label{estim v cas 1}
\left| v_{j,k}^{(i)} \right| \le C_\varepsilon e^{-k^2T}\| y^0 \|_{L^2(0,\pi)^2}, 
 \quad i,j=1,2.
\end{equation}

\smallskip

\textbf{Case 2: }Let $k\in \Lambda_{1}$ such that $k > k_\varepsilon$ 
and $ \left| I_k (p,q) \right|^{-1} \le e^{k^2(T_0 (p,q) + 2 \varepsilon)}$. \\
As in the previous case, we take $v^{(2)}_{1,k}=v^{(2)}_{2,k}=0$  
and the moment problem \eqref{moment Lambda1} has a unique solution, 
given by \eqref{sol moment case 1}. 
Thanks to the property of $\psi_k^*$ (see \eqref{estm psi2}) and Lemma \ref{lemme fk}, we get for $i=1,2$ the following estimates
\begin{equation}\label{estim Ik psik}
|f^{(1)}_k|\geqslant C/k^3,~~~
|\widetilde{f}^{(i)}_k|\leqslant \dfrac{C}{k},~~~
|y_{i,k}^0|\leqslant C\|y^0\|_{L^2(0,\pi)^2},~~~\forall k\in\mathbb{N}^*.
\end{equation}
Thus, using the assumptions on $k$, we obtain
\begin{equation*}
\left\{\begin{array}{l}
|v_{1,k}^{(1)}|\le C k^3e^{- k^2 T} \| y^0 \|_{L^2(0,\pi)^2} 
\le C_\varepsilon e^{-(T- \varepsilon) k^2}\| y^0 \|_{L^2(0,\pi)^2},\\  \noalign{\smallskip}
|v^{(2)}_{1,k}|\le
 \dfrac{C_\varepsilon e^{-(T- \varepsilon) k^2}}{|I_k(p,q)|}\| y^0 \|_{L^2(0,\pi)^2}\le
C_\varepsilon e^{-(T-T_0- 3\varepsilon) k^2}\| y^0 \|_{L^2(0,\pi)^2},
\end{array}\right.
\end{equation*}
where $C_\varepsilon$ is a constant which is independent on $k$ and $y^0$.

\textbf{Case 3: }Consider now $k\in \Lambda_{1}$ such that $k > k_\varepsilon$ 
and $ \left| I_k (p,q) \right|^{-1} > e^{k^2(T_0 (p,q) + 2 \varepsilon)}$. \\
This implies with \eqref{kepsilon} that 
\begin{equation}\label{esti Ik res}
\left| I_{a,k} (p,q) \right|^{-1} < e^{k^2(T_0 (p,q) +  \varepsilon)}.  
\end{equation}
The two last inequalities  lead to
\begin{equation*}
\left|I_{k}(p,q)\right|<e^{-\varepsilon k^2}
\left|I_{a,k}(p,q)\right|.
\end{equation*}
Combined with inequality \eqref{cond prop intersec2}, taking $k_{\varepsilon}$ large enough, we get
\begin{equation}\label{detAk}
\displaystyle \left| \det A_{1,k} \right| 
> C_\varepsilon e^{- \varepsilon k^2} \left|I_{a,k} ( p,q )\right|,
\end{equation}
with $C_\varepsilon$ independent on $k$. To solve the moment problem \eqref{moment Lambda1}, 
we take here $v_{2,k}^{(1)}=v_{2,k}^{(2)}=0$. 
Then the moment problem \eqref{moment Lambda1} reads $A_{1,k} V_{1,k} = F_k$. 
Since $\mathrm{det}A_{1,k}\neq0$ and using \eqref{detAk}, the inverse of $A_{1,k}$ is given by
\begin{equation*}
(A_{1,k})^{-1} = (\mathrm{det}~A_{1,k})^{-1}\left(
\begin{array}{cc}
  f_k^{(2)}& - \widetilde f_k^{(2)}\\ 
-f_k^{(1)} & \widetilde f_k^{(1)}
\end{array}
\right).
\end{equation*}
We deduce that the solution to the moment problem \eqref{moment Lambda1} is
\begin{equation*}
\left\{\begin{array}{l}
v_{1,k}^{(1)}=\frac{e^{- k^2 T}}{\mathrm{det}~A_{1,k}}
\{-f_k^{(2)}y_{1,k}^0
+(TI_k(p,q)f_k^{(2)}+\widetilde{f}_{k}^{(2)})y_{2,k}^0\},\\  \noalign{\smallskip}
v^{(2)}_{1,k}=
\frac{e^{- k^2 T}}{\mathrm{det}~A_{1,k}}
\{f_k^{(1)}y_{1,k}^0
-(TI_k(p,q)f_k^{(1)}+\widetilde{f}_{k}^{(1)})y_{2,k}^0\}.
\end{array}\right.
\end{equation*}
The last expression together with \eqref{esti Ik res} and \eqref{detAk}    gives
\begin{equation}\label{resol moment cas 3}
|v_{1,k}^{(i)}|
\le C_\varepsilon e^{-(T-T_0-2\varepsilon) k^2}\| y^0 \|_{L^2(0,\pi)^2}, \quad i=1,2.
\end{equation}

\textbf{Case 4:} Let  us consider $k \in \Lambda_{2}$.\\
If $k\le k_\varepsilon$, we can argue as in Case 1. Let us suppose that $k>k_\varepsilon$. 
In this case, $I_{a,k}(p,q)=0$, $I_k(p,q)\neq0$ and inequality 
\eqref{kepsilon} reads 
$\left|I_{k} ( p,q ) \right|^{-1} <e^{ k^2(T_0 (p,q) + \varepsilon)}$. 
We  take here $v^{(2)}_{1,k}=v^{(2)}_{2,k}=0$ 
and the solution of moment problem \eqref{moment Lambda1} is given by \eqref{sol moment case 1}. 
We get
\begin{equation*}
\left| v_{j,k}^{(i)} \right| \le C_\varepsilon e^{-k^2(T-T_0(p,q)-2\varepsilon)}\| y^0 \|_{L^2(0,\pi)^2}, 
\quad  i,j=1,2.
\end{equation*}

 \textbf{Case 5:} Let us  now deal with the case $k \in \Lambda_3$.   \\
 We recall that $I_k(p,q)=0$, $I_{1,k}(p,q)\neq0$ and inequality 
\eqref{kepsilon} reads 
\begin{equation}\label{resol moment:cas5}
\left|I_{a,k} ( p,q ) \right|^{-1} <e^{ k^2(T_0 (p,q) + \varepsilon)}.
\end{equation} 
  The moment problem \eqref{moment Lambda1} 
 is now $A_{1,k}V_{1,k}=F_k$ with $A_{1,k}$ and $F_k$ given in \eqref{Ak} and \eqref{Fk}, 
 respectively. From  \eqref{cond prop intersec2}, the matrix $A_{1,k}$ is invertible and
 \begin{equation*}
\left\{\begin{array}{l}
v_{1,k}^{(1)}=\frac{e^{- k^2 T}}{\mathrm{det}~A_{1,k}}
\{-f_k^{(2)}y_{1,k}^0
+\widetilde{f}_{k}^{(2)}y_{2,k}^0\},\\  \noalign{\smallskip}
v^{(2)}_{1,k}=\frac{e^{- k^2 T}}{\mathrm{det}~A_{1,k}}
\{f_k^{(1)}y_{1,k}^0-\widetilde{f}_{k}^{(1)}y_{2,k}^0\}.
\end{array}\right.
\end{equation*}
Using inequalities \eqref{cond prop intersec2} and \eqref{resol moment:cas5}, we obtain estimate \eqref{resol moment cas 3}.

\textbf{Conclusion:} \\
We have constructed a control  $v$ of the form  \eqref{contv} and \eqref{contvi}, 
which satisfies
\begin{equation*}
\left| v_{j,k}^{(i)} \right| \le C_\varepsilon 
e^{-k^2(T-T_0(p,q)-3\varepsilon)}\| y^0 \|_{L^2(0,\pi)^2}, 
\quad i,j=1,2,~k\in\mathbb{N}^*.
\end{equation*}
The last inequality, the estimate \eqref{estim qik} of $q_{i,k}$ and 
the expression  \eqref{contvi} of $v^{(i)}$ ($i=1,2$) lead to
\begin{equation*}
\left\| v^{(i)} \right\|_{L^2(0, T)} 
\le C_{\varepsilon ,T}  e^{ -k^2 ( T - T_0 (p,q) - 4 \varepsilon)},~~i=1,2.
\end{equation*}
Thus, taking 
$\varepsilon \in (0, (T - T_0 (p,q))/4)$, 
we have the absolute convergence of the series defining $v^{(1)}$ and $v^{(2)}$ in $L^2 (0, T)$. 
This ends the proof.

\end{proof}

\begin{proof}[Proof of Theorem \ref{theo null int}]
Using Proposition \ref{prop intersec2}, System \eqref{system primmal int} is equivalent to a system with coupling terms $\widehat{p}$ and $\widehat{q}$ satisfying Condition 
\eqref{cond approx2} and \eqref{cond prop intersec2}. 
Proposition \ref{prop resol} leads to the null controllability of System \eqref{system primmal int} when $T>T_0(\widehat{p},\widehat{q})$. 
We end the proof of Theorems \ref{theo null int} remarking that $T_0(\widehat{p},\widehat{q})=0$.

\end{proof}

\section{Proof of Theorem \ref{theo null T0}}

\subsection{Positive null controllability result}


 Before  studying the case where the intersection of the coupling and control domains 
is empty, 
we will first rewrite the function $\psi_k^*$ 
given in  Proposition \ref{prop base}. 

\begin{lmm}\label{lemme psi_k} Let $k\in\mathbb{N}^*$. 
 Consider the function $\psi_k^*$ defined in Proposition \ref{prop base}. 
 If we suppose that Condition \eqref{cond p q non nuls} holds, 
 then for all $x\in\omega$
 \begin{equation*}
  \psi_k^*(x)=\tau_k\varphi_k(x)+g_k(x)\mathrm{~for~all}~x\in\omega,
 \end{equation*}
where 
\begin{equation*}
\left\{\begin{array}{l} \tau_k:=   \alpha_k^*-\sqrt{\dfrac{\pi}{2}}\dfrac{1}{k}\displaystyle\int_0^a\cos(k\xi)
       [\partial_x(p(\xi)\varphi_k(\xi))-q(\xi)\varphi_k(\xi)]d\xi,\\
 g_k(x):=   -\dfrac {I_k(p,q)}{k}\displaystyle\int_0^x\sin(k(x-\xi))\varphi_k(\xi)d\xi     
 -\sqrt{\dfrac{\pi}{2}}\dfrac{1}{k}I_{a,k}(p,q)\cos(kx).
\end{array}\right.\end{equation*}
\end{lmm}

\begin{proof}
Since $p=q\equiv0$ in $\omega$,  we get for all $x\in \omega$,
\begin{equation*}
\begin{array}{rcl}
\displaystyle\psi_k^* (x) &=& \alpha_k^* \varphi_k (x) 
- \dfrac {I_k(p,q)}{k} \displaystyle\int_0^x \sin (k (x- \xi)) \varphi_k (\xi) \, d \xi \\
&&\hspace*{4cm}
- \dfrac 1k \displaystyle\int_0^a\sin (k (x-\xi))
[\partial_x(p(\xi)\varphi_k(\xi))-q(\xi)\varphi_k (\xi) ]d\xi. \\
\end{array}
\end{equation*}

\end{proof}

\begin{proof}[Proof of Theorem \ref{theo null T0}]
We will follow the strategy of \cite{Ammar-Khodja2015}. 
More precisely, we will prove Theorem 1.3 with the help of Proposition \ref{prop resol}.
Assume that Conditions \eqref{cond approx2} and \eqref{cond p q non nuls} hold. 
Consider the functions $ f^{(1)}$ and $ f^{(2)}$ defined in Lemma \ref{lemme fk} 
and the matrix $A_{1,k}$ given in \eqref{Ak}. 
Let $k\in\mathbb{N}^*$. 
We recall that
\begin{equation*}
\det A_{1,k} =  \widetilde{f}_{k}^{(1)}  f_k^{(2)}-\widetilde  f_k^{(2)} f_k^{(1)}, 
\end{equation*}
where, for $i=1,2$, $f_{k}^{(i)}$ and $\widetilde{f}_{k}^{(i)}$ are defined in \eqref{def fik}. 
Since $\mathrm{Supp}~f^{(i)}\subseteq\omega$, using the expression of $\psi_k^*$ given in Lemma \ref{lemme psi_k}, 
we obtain 
\begin{equation*}
\widetilde f_{k}^{(i)}=\tau_kf_{k}^{(i)}+\int_0^\pi f^{(i)}(x)\,g_k(x)\,dx,
\end{equation*}
where for all $x\in \omega$
\begin{equation*}
\displaystyle g_k (x)=-\frac{I_k(p,q)}{k}\int_{0}^{x}\sin( k(x-\xi ) ) \varphi _{k} ( \xi  )d\xi
-\sqrt{\frac\pi 2}\frac{1}{k}  I_{a,k}(p,q)  \cos(kx) .
 \end{equation*}
We deduce that
\begin{equation*}\begin{array}{rcl}
\det A_{1,k}&=&f_k^{(2)}\displaystyle\int_0^\pi f^{(1)}(x)\,g_k(x)dx
-f_k^{(1)}\displaystyle\int_0^\pi f^{(2)}(x)\,g_k(x)\,dx \\ \noalign{\smallskip}
&=&- \dfrac {I_k(p,q)} {k} \  \left( 
f_k^{(2)} \displaystyle\int_0^\pi \int_{0}^{x}  f^{(1)} (x) \sin(k(x-\xi))\varphi _{k} (\xi ) d\xi\,dx\right.\\  \noalign{\smallskip} 
&&\hspace*{0.5cm}\left.-f_k^{(1)}\displaystyle\int_0^\pi \int_{0}^{x}  f^{(2)} (x)\sin(k(x-\xi))\varphi _{k} (\xi)d\xi\,dx \right)- \sqrt{\dfrac \pi 2}\dfrac 1k  I_{a,k}(p,q) \left(\widehat f_k^{(1)}{f}_{k}^{(2)}-\widehat f_k^{(2)}{f}_{k}^{(1)} \right),
\end{array} \end{equation*}
where $\widehat f_{k}^{(i)}$ are defined in \eqref{def hat fik}. 
 Since the integrals 
\begin{equation*}
\int_0^\pi \int_{0}^{x}  f^{(i)} (x) \sin  ( k(x-\xi ) ) \varphi _{k} ( \xi  ) \,d\xi\,dx
\end{equation*}
and the sequence $(f^{(i)}_k)_{k\in\mathbb{N}^*,i\in\{1,2\}}$ are uniformly bounded with respect to $k$ and $i$, we conclude with the help of Lemma \ref{lemme fk}.

We deduce that Condition \eqref{cond prop intersec2} holds. Thus, using Proposition \ref{prop resol}, 
 System \eqref{system primmal int} is null controllable at time $T$.

\end{proof}

\subsection{Negative null controllability result}

 Let us now prove the negative part of Theorem 1.3 
with the  strategy used in \cite{Ammar-Khodja2015}. Let $T<T_0(p,q)$. 
We will argue by contradiction: 
Assume that System \eqref{system primmal int} is null controllable at time $T$. 
Using Proposition \ref{prop ine obs}, there exists a constant $C_{obs}>0$ such that 
for all $\theta^0\in L^2(0,\pi)^2$, the solution to System \eqref{system dual} satisfies 
the observability inequality 
 \begin{equation}\label{null cont neg}
 \|\theta(0)\|^2_{L^2(0,\pi)^2}\leqslant C_{obs} \iint_{Q_T}|\mathds{1}_{\omega}(x)B^*\theta(x,t)|^2dxdt.
 \end{equation}
Using the Definition  of $T_{0}(p,q)$ (see \eqref{def T_0}) 
there exists a strictly increasing  sequence   
$(k_n)_{n\in\mathbb{N}^*} \subseteq \mathbb{N}$ satisfying: 
\begin{equation}\label{negative T_0}
T_0 (p,q)=\lim_{n\to\infty} 
\frac{\min \left(\log\left|I_{a,k_n}(p,q)^{-1}\right|,\log\left|I_{k_n}(p,q)^{-1}\right|\right)}{k_n^2}. 
\end{equation}
Let us fix $n\geq 1$ and  $\theta^{0n} := a_{n} \Phi_{1,k_n}^{\ast } + b_{n}\Phi _{2,k_n}^{\ast }$ 
with  $(a_{n},b_{n})\in \mathbb{R}^{2}$ to be determined later and $\Phi_{2,k_n}^{\ast }$, $\Phi_{1,k_n}^{\ast }$ the eigenfunction 
and generalized eigenfunction associated with $k_n^2$ 
given in Proposition~\ref{prop base}. 
If we denote by $\theta^n$ the solution to the dual System \eqref{system dual} for  initial data 
$\theta^{0n}$, then
\begin{equation*}
 \theta^n(x,t)=e^{-k_n^2(T-t)}\{a_n\Phi^*_{1,k_n}+(b_n-(T-t)I_{k_n}(p,q)a_n)\Phi^*_{2,k_n}\},
\end{equation*}
thus, using the orthogonality $\langle\psi_{k_n}^*,\varphi_{k_n}\rangle_{L^2(0,\pi)}=0$, we have
\begin{equation*}
\left\{\begin{array}{l}
D_{1,n}:=\|\theta^n(0)\|_{L^2(0,\pi)^2}^2=\displaystyle e^{-2k_n^{2}T} \left\{ |a_{n}|^{2} |\psi _{k_n}^*|^{2}
  + \left( b_n -TI_{k_n}(p,q) a_n \right)^2 + |a_{n}|^{2}  \right\} \\
D_{2,n}:= \iint_{Q_T}|\mathds{1}_{\omega}(x)B^*\theta^n(x,t)|^2dxdt= 
  \displaystyle\int_{0}^T\displaystyle\int_{\omega}e^{-2k_n^{2}t}
 \left| a_{n} \psi _{k_n}^*(x)+(b_{n}- tI_{k_n}(p,q)a_{n})
 \varphi_{k_n}(x)\right|^{2}\,dx \,dt.
\end{array}\right.\end{equation*}
The observability inequality \eqref{null cont neg} reads
\begin{equation}\label{ine obs neg2}
 D_{1,n}\leqslant C_{obs}D_{2,n}.
\end{equation}
By choosing  $a_n: = 1$ and $b_{n}:=- \tau _{k_n}$, we get
\begin{equation}\label{est A1n}
\begin{array}{rcl}
D_{1,n}&\geqslant& e^{-2k_n^2T}
\end{array}\end{equation}
and the expression of  $\psi_{k_n}^*(x)$ given in Lemma \ref{lemme psi_k} 
leads to
\begin{equation*}
\begin{array}{rcl}
D_{2,n}&=&\displaystyle\int_{0}^T\displaystyle\int_{\omega} e^{-2k_n^{2}t} \left|
- \sqrt{\frac \pi 2}  \frac{1}{k_n}I_{a,k_n}(p,q) \cos (k_nx)\right.\\\noalign{\smallskip}
&&\hspace*{2.5cm}\left.- I_{k_n}(p,q)\dfrac{1}{k_n}\displaystyle \int_0^x 
\sin \left(k_n (x- \xi) \right) \varphi_{k_n} (\xi) \, d\xi
-tI_{k_n}(p,q) \varphi _{k_n}(x) \right|^{2} \,dx \,dt \\\noalign{\smallskip}
&\leqslant& C(I_{a,k_n}(p,q)^2+I_{k_n}(p,q)^2).
\end{array}\end{equation*}
Let $\varepsilon >0$. Equality \eqref{negative T_0} implies that 
there is $k_{\varepsilon}\in\mathbb{N}^*$ such that 
for all $k_n \geqslant k_{\varepsilon}$
\begin{equation*}
\max \left(\left|I_{a,k_n}(p,q)\right|^{2},\left|I_{k_n}(p,q)\right|^{2}\right) 
\leqslant e^{-2k_n^2(T_0(p,q)-\varepsilon)}.
\end{equation*}
We deduce that for $\varepsilon:=(T_0(p,q)-T)/2$, we get
\begin{equation}\label{est A2n}
 D_{2,n}\leqslant Ce^{-2k_n^2(T+\varepsilon)}.
\end{equation}
Thus, since $k_n$ goes to $\infty$, 
  estimates \eqref{est A1n} and \eqref{est A2n} are in contradiction with 
 inequality \eqref{ine obs neg2} for $n$ large enough.

\section{Proof of Theorem 1.2}

We will proved Theorem \ref{theo approx int} using the criterion of Fattorini, 
as in the pioneer work \cite{olive_bound_appr_2014}.

\begin{thrm}[see \cite{fattorini66}, Cor. 3.3]\label{theo hautus}
 System \eqref{system primmal int} is approximatively controllable at time $T$ if and only if 
 for any $s\in\mathbb{C}$ and for any $u\in \mathcal{D}(L^*)$ we have
 \begin{equation*}
\left.\begin{array}{lll}  L^*u=su&\mathrm{~in~}&(0,\pi)\\ \noalign{\smallskip}
       B^*u=0&\mathrm{~in~}&\omega\end{array}\right\}\Rightarrow u=0.
 \end{equation*}

\end{thrm}

\begin{proof}[Proof of Theorem \ref{theo approx int}]$\left.\right.$

\textbf{Necessary condition: }
Let us suppose that Conditions \eqref{cond approx1}-\eqref{cond approx2} do not hold
 \emph{i.e.} there exists $k_0\in\mathbb{N}^*$ such that 
 \begin{equation*}\begin{array}{c}
  I_{k_0}(p,q)=I_{a,k_0}(p,q)=0
  \mathrm{~and~}
   (\mathrm{Supp}~p\cup\mathrm{Supp}~q)\cap\omega=\varnothing.
 \end{array}\end{equation*}
We remark that   the function $\psi_{k_0}^*$ of Lemma \ref{lemme psi_k} 
satisfy $\psi_{k_0}^*=\tau_{k_0}\varphi_{k_0}$ in $\omega$, then
$$\Phi^*_{1,k_0}-\tau_{k_0}\Phi_{2,k_0}^*=\left(\begin{array}{c}0\\\varphi_k\end{array}\right)
\mbox{ in }\omega.$$ 
We deduce that  $\Phi^*_{1,k_0}-\tau_{k_0}\Phi_{2,k_0}^*$ 
is an non-trivial eigenfunction associated with the eigenvalue $k_0^2$ 
 of the operator $L^*$ satisfying
 \begin{equation*}
  B^*(\Phi^*_{1,k_0}-\tau_{k_0}\Phi_{2,k_0}^*)\equiv0 \mathrm{~in~}\omega.
 \end{equation*}
Thus, using Theorem \ref{theo hautus}, 
System \eqref{system primmal int} is not 
  approximately controllable at time $T$.

\textbf{Sufficient condition: }
Let us suppose that Conditions \eqref{cond approx1}-\eqref{cond approx2}  hold.
If $(\mathrm{Supp}~p\cup\mathrm{Supp}~q)\cap\omega\neq\varnothing$, then 
we conclude using Theorem \ref{theo null int}. Let us now suppose that 
 \begin{equation*}
   (\mathrm{Supp}~p\cup\mathrm{Supp}~q)\cap\omega=\varnothing
     \hspace*{5mm} \mbox{and}\hspace*{5mm} 
    |I_{k}(p,q)|+|I_{a,k}(p,q)|\neq0\mathrm{~for~all~}k\in \mathbb{N}^*.
\end{equation*}
 
If $I_k(p,q)\neq0$, the set of the eigenvectors associated with the eigenvalue $k^2$ of $L^*$ is generated  
by $\Phi_{2,k}^*$ (see Proposition \ref{prop base}). 
In this case, we remark that for all $k\in \mathbb{N}^*$
\begin{equation}\label{preuve approx}
 B^*\Phi_{2,k}^*=\varphi_k\not\equiv0\mathrm{~in~}\omega.
\end{equation}

 If $I_k(p,q)=0$, the eigenvectors associated with the eigenvalue $k^2$ of $L^*$ are linear combinations of 
 $\Phi_{1,k}^*$ and $\Phi_{2,k}^*$. 
 Let $\alpha,\beta\in\mathbb{R}$ and $\Phi^*:=\alpha\Phi_{1,k}^*+\beta\Phi_{2,k}^*$ satisfying
\begin{equation}\label{preuve approx bis}
  B^*\Phi^*\equiv0\mathrm{~in~}\omega.
\end{equation} 
 Using Lemma \ref{lemme psi_k}, it is equivalent to
 \begin{equation*}
 (\alpha+ \beta\tau_k)\varphi_k(x)-\beta\sqrt{\frac{\pi}{2}}\frac{1}{k}I_{a,k}(p,q)\cos(kx)=0\mathrm{~for~all~}x\in\omega.
 \end{equation*}
 Since $I_{a,k}(p,q)\neq0$, we deduce that $\beta=0$. Then $\alpha=0$. 
 We conclude with the help of Theorem \ref{theo hautus}.

\end{proof}

\section{Proof of Theorem \ref{theo bord}}

 As in Section \ref{moment problem}, System \eqref{system primmal bord} 
is null controllable at time $T$ if and only if for all 
$y^0\in H^{-1}(0,\pi)^2$, $k\in \mathbb{N}^*$ and $i\in\{1,2\}$ the solution 
$\theta_{i,k}$ to the dual System \eqref{system dual} for the initial data $\Phi_{i,k}^*$ satisfies
\begin{equation}\label{equiv null bord}
 \displaystyle\int_0^Tu(t)B^*\partial_x\theta_{i,k}(0,t)dt
 =-\langle y^0,\theta_{i,k}(\cdot,0)\rangle_{H^{-1},H^1_0}.
\end{equation}
We recall that, for all $k\in \mathbb{N}^*$, 
$\theta_{1,k}$ and $\theta_{2,k}$ are given for all $(x,t)\in Q_T$ by
\begin{equation*}
\theta_{1,k}(x,t)=e^{-k^2 (T-t)}\left(\Phi_{1,k}^*(x)-(T-t)I_k(p,q)\Phi_{2,k}^*(x)\right)\mbox{~and~}
\theta_{2,k} (x, t) = e^{-k^2 (T-t)}\Phi_{2,k}^*(x). 
\end{equation*}

\begin{proof}[Proof of Theorem \ref{theo bord}]
Again, we will follow the strategy used in \cite{Ammar-Khodja2015}. 
 Assume that $T>T_1$ and $I_k(p,q)\neq 0$ for all $k\in\mathbb{N}^*$. 
 We will look for the control $u$ under the form 
\begin{equation}\label{expr cont bord}
u(t):=\sum\limits_{k\in\mathbb{N}^*}\{u_{1,k}q_{1,k}(T-t)+u_{2,k}q_{2,k}(T-t)\},
\end{equation}
for all $t\in(0,T)$, where $q_{1,k}$ and $q_{2,k}$ are defined in Section \ref{moment problem}. 
Plugging the expressions of $u$, $\theta_{1,k}$ and $\theta_{2,k}$ 
in Equality \eqref{equiv null bord}, we obtain the \textit{moment problem} 
\begin{equation*}
 \left\{\begin{array}{l}
u_{1,k}=-e^{-k^2T}\dfrac{\langle y_1^0,\varphi_k\rangle_{H^{-1},H^1_0}}{\partial_x\varphi_k(0)},\\  \noalign{\smallskip}
u_{2,k}=\dfrac{e^{-k^2T}}{I_k\partial_x\varphi_k(0)}
\{\langle y_1^0,\psi_k^*\rangle_{H^{-1},H^1_0}
+\langle y_2^0,\varphi_k\rangle_{H^{-1},H^1_0}
-\left(I_kT+\dfrac{\partial_x\psi_k^*(0)}{\partial_x\varphi_k(0)}\right)
\langle y_1^0,\varphi_k\rangle_{H^{-1},H^1_0}   \} .
        \end{array}
\right.
\end{equation*}
Let $\varepsilon>0$. Using  the definition  of $T_1$ (see \eqref{def T1234}), 
we have $I_k(p,q)>C_{\varepsilon}e^{-k^2(T_1+\varepsilon)}$ for all $k\in \mathbb{N}^*$. Then, using the estimates \eqref{estm psi2} and \eqref{estim Ik psik}, we get
\begin{equation*}
 |u_{1,k}|+|u_{2,k}|\leqslant Ce^{-k^2(T-T_1-2\varepsilon)}\|y^0\|_{H^{-1}(0,\pi)^2}.
\end{equation*}
Thus for $\varepsilon < (T-T_1)/2$, the control $u$ defined in \eqref{expr cont bord} is an element 
of $L^2(0,T)$.

 Assume now that $T<T_1$ and $I_k(p,q)\neq 0$ for all $k\in\mathbb{N}^*$. By contradiction let us suppose that 
 there exists a constant $C_{obs}$ such that for all $\theta^0\in H^{1}_0(0,\pi)^2$ the solution to the dual System 
 \eqref{system dual} satisfies
  \begin{equation}\label{null cont bord 2}
 \|\theta(0)\|^2_{H^1_0(0,\pi)^2}\leqslant C_{obs} \int_{0}^T|B^*\partial_x\theta(0,t)|^2dt.
 \end{equation} 
 Let $\varepsilon=(T_1-T)/2$. Using the definition of $T_1$, there exists a sequence $(k_n)_{n\in\mathbb{N}^*}$ such that 
 \begin{equation}\label{ine Ik bord}
  I_{k_n}(p,q)<e^{-k^2_n(T+\varepsilon)}.
 \end{equation}
 Let $\theta^0_n:=a_n\Phi_{1,k_n}^*+b_n\Phi_{2,k_n}^*$ with $(a_n,b_n)\in\mathbb{R}^2$. We recall that
 \begin{equation*}
 \theta^n(x,t)=e^{-k_n^2(T-t)}\{a_n\Phi^*_{1,k_n}+(b_n-(T-t)I_k(p,q)a_n)\Phi^*_{2,k_n}\}.
 \end{equation*}
Then,  after calculation, we get
 \begin{equation*}
   \|\theta(0)\|^2_{H^1_0(0,\pi)^2}=e^{-2k_n^2T}
   (a_n^2\|\psi_{k_n}\|_{H^1_0}^2+a_n^2k_n^2+(b_n-TI_k(p,q)a_n)^2k_n^2)
 \end{equation*}
and
  \begin{equation*}
   \int_{0}^T|B^*\partial_x\theta(0,t)|^2dt
   =  \int_{0}^Te^{-2k_n^2(T-t)}|a_n\partial_x\psi_{k_n}(0)+\sqrt{\frac{2}{\pi}}(b_n-(T-t)I_{k_n}(p,q)a_n)k_n|^2dt.
 \end{equation*}
For $a_n:=1$ and $b_n:=-\sqrt{\frac{\pi}{2}}\partial_x\psi_{k_n}(0)/k_n$, taking into account  inequality \eqref{ine Ik bord} and using the estimate \eqref{estm psi2}, 
we obtain
\begin{equation*}
    \|\theta(0)\|^2_{H^1_0(0,\pi)^2}\geqslant k_n^2e^{-2k_n^2T}
  \hspace*{5mm}   \mathrm{and}\hspace*{5mm} 
       \int_{0}^T|B^*\partial_x\theta(0,t)|^2dt\leqslant C k_n^2e^{-2k_n^2(T+\varepsilon)}.
\end{equation*}
Thus for $n$ large enough we get a contradiction with observability inequality \eqref{null cont bord 2}.

 \end{proof}

\section{Comments and open problems}

When the control domain and the support of the coupling coefficients $p$ and $q$ is disjoint in the system 
\begin{equation}\label{system primmal concl}
 \left\{\begin{array}{ll}
\partial_ty_1-\partial_{xx}y_1=\mathds{1}_{\omega}v
&\mathrm{in~} Q_T,\\
\partial_ty_2-\partial_{xx}y_2+p(x)\partial_xy_1+q(x)y_1=0
&\mathrm{in~} Q_T,\\
\noalign{\smallskip}y_1(0,\cdot)=y_1(\pi,\cdot)=y_2(0,\cdot)=y_2(\pi,\cdot)=0&\mathrm{on}~(0,T),\\
\noalign{\smallskip}y_1(\cdot,0)=y^0_1,~y_2(\cdot,0)=y^0_2&\mathrm{in}~(0,\pi)
        \end{array}
\right.
\end{equation}
(resp. system \eqref{system primmal bord}), it is legitimate to ask if the minimal time $T_1$ (resp. $T_0$)  given  in Theorem 1.3 (resp. Theorem 1.4) 
can be different of zero and finite. For $p\equiv0$ in $(0,\pi)$, it is proved in 
 \cite[Lem. 7.1]{Ammar-Khodja2015} that for any $\tau_0\in[0,\infty]$ there exists a function $q\in L^{\infty}(0,\pi)$ 
such that the minimal time of null controllability $T_0(p,q)$ 
associated with System  \eqref{system primmal int} is given by $T_0(p,q)=\tau_0$.  
The authors give explicit functions and one can  easily adapt them to the case $p\not\equiv0$ in $(0,\pi)$. 
In the other hand, the null controllability in the cases $T=T_0$ in Theorem \ref{theo null T0} and $T=T_1$ in Theorem \ref{theo bord} are  open problems.

In higher space dimension, even for this simplified system \eqref{system primmal concl} (resp. system \eqref{system primmal bord}), 
distributed and boundary controllability are also open problems. 
Considering the different results described in the introduction of the present paper, we can conjecture that  the  system of two coupled linear parabolic equations
\begin{equation}\label{system primmal gene concl}
 \left\{\begin{array}{ll}
\partial_ty_1=\Delta y_1+g_{11}\cdot\nabla y_1+g_{12}\cdot\nabla y_2 +a_{11}y_1+a_{12}y_2+\mathds{1}_{\omega}v
&\mathrm{in~} \Omega\times(0,T),\\
\noalign{\smallskip}\partial_ty_2=\Delta y_2+g_{21}\cdot\nabla y_1+g_{22}\cdot\nabla y_2 +a_{21}y_1+a_{22}y_2
&\mathrm{in~} \Omega\times(0,T),\\
\noalign{\smallskip}y=0&\mathrm{on}~\partial\Omega\times(0,T),\\
\noalign{\smallskip}y(\cdot,0)=y^0&\mathrm{in}~\Omega,
        \end{array}
\right.
\end{equation}
is null controllable at time $T>0$ if there exists an open nonempty subset $\omega_0$ of $\omega$ such that
\begin{equation}\label{concl:conj}
|a_{21}|>C\mbox{ in }\omega_0\times(0,T)
\hspace*{5mm} \mbox{or}\hspace*{5mm} 
|g_{21}^k|>C\mbox{ in }\omega_0\times(0,T),
\end{equation}
for a $k\in\{1,...,N\}$. 

It seems that the main difficulty is to prove a Carleman estimate for the adjoint problem of system \eqref{system primmal gene concl} under condition \eqref{concl:conj} when the coupling term is a differential operator (see for instance  \cite{benabdallah2014,guerrerosyst22} and also \cite{pierreordre1} for a different approach). 
In the one-dimensional case, we were not able to adapt the strategy developed in this paper 
in this general setting.

\begin{acknowledgement} The author thanks  Assia Benabdallah, Manuel González-Burgos 
and Farid Ammar Khodja   for their interesting comments and suggestions. 
He thanks as well  the  referees for his remarks that helped to improve the paper.
\end{acknowledgement}


\end{document}